\documentclass[12pt]{article}

\usepackage{todonotes}
\presetkeys{todonotes}{inline}{}

\usepackage{amsmath,amssymb,latexsym, amsbsy, amsfonts, amscd, amsthm, mathrsfs, xy,comment, xcolor}
\usepackage{multirow}
\usepackage{hyperref}
\usepackage{enumitem, empheq, url}
 \usepackage[bibtex-style]{amsrefs}
\xyoption{all}
\usepackage{tikz-cd}
\usepackage{rotating, scalerel}

\theoremstyle{plain}

\newtheorem{theorem}{Theorem}[section]

\newtheorem{prop}[theorem]{Proposition}
\newtheorem{corollary}[theorem]{Corollary}

\newtheorem{lemma}[theorem]{Lemma}

\theoremstyle{definition}
\newtheorem{definition}[theorem]{Definition}
\newtheorem{remark}[theorem]{Remark}
\newtheorem{example}[theorem]{Example}

\long\def\symbolfootnote[#1]#2{\begingroup
\def\thefootnote{\fnsymbol{footnote}}\footnote[#1]{#2}\endgroup}

\def\lra{\longrightarrow}

\DeclareMathOperator{\GL}{GL}

\def\D{\mathbf{D}}

\def\cG{\mathcal{G}}
\def\cI{\mathcal{I}}
\def\cA{\mathcal{A}}

\def\sgn{\mathrm{sgn}}
\def\N{\mathrm{N}}
\def\1{\mf{1}}

\DeclareMathOperator{\Fitt}{Fitt}

\DeclareMathOperator{\Sym}{Sym}
\DeclareMathOperator{\Ind}{Ind}

\DeclareMathOperator{\Frac}{Frac}

\DeclareMathOperator{\Hom}{Hom} 

\DeclareMathOperator{\sign}{sign}
 
\DeclareMathOperator{\ord}{ord}

 \DeclareMathOperator{\real}{Re}

 \DeclareMathOperator{\red}{red}

\DeclareMathOperator{\Spec}{Spec}
\DeclareMathOperator{\Ext}{Ext}

\DeclareMathOperator{\tr}{tr}

\DeclareMathOperator{\im}{im}

\newcommand{\mat}[4]{\begin{pmatrix} {#1} & {#2} \\ {#3} & {#4} \end{pmatrix}}

\newcommand{\mf}{\mathfrak }

\def\fp{\mathfrak{p}}
\def\fq{\mathfrak{q}}

\def\fm{\mathfrak{m}}

\def\fP{\mathfrak{P}}

\def\fm{\mathfrak{m}}

\def\T{\mathbf{T}}
\def\Z{\mathbf{Z}}
\def\F{\mathbf{F}}
\def\Q{\mathbf{Q}}
\def\C{\mathbf{C}}
\def\G{\mathbf{G}}
\def\B{\mathbf{B}}

\def\bdf{\begin{defn}}
\def\edf{\end{defn}}

\def\cO{\mathcal{O}}

\def\cP{\mathcal{P}}

\def\cA{\mathcal{A}}

\def\cV{\mathcal{V}}

\def\cG{\mathcal{G}}

\def\fb{\mathfrak{b}}

\def\Gal{{\rm Gal}}

\def\cB{{\cal B}}

\def\cP{{\cal P}}

\def\ram{\text{ram}}

\begin{document}
\baselineskip 15.8pt

\title{The Residually Indistinguishable Case of Ribet's Method for $\GL_2$}
\author{Samit Dasgupta \\ Mahesh Kakde \\ Jesse Silliman \\ Jiuya Wang}
\maketitle
\tableofcontents

\bigskip

\section{Introduction}

Since the publication of the groundbreaking paper \cite{ribet}, Ribet's method has played a central role in modern algebraic number theory.
Ribet's method provides a strategy for constructing nontrivial extensions of a $p$-adic Galois representation $\rho_1$ by another such representation $\rho_2$.  After Ribet's original proof of the converse to Herbrand's theorem, his method was used by Mazur--Wiles to prove the main conjecture of Iwasawa theory over $\Q$ \cite{mw} and by Wiles to prove the main conjecture over arbitrary totally real fields \cite{wiles}.  An important innovation introduced in Mazur--Wiles is the use of Fitting ideals. More recently, generalizations of Ribet's method have been used by Skinner--Urban to prove the main conjecture for $p$-ordinary elliptic curves \cite{su} and by the first two authors of this paper to prove the Brumer--Stark conjecture away from $p=2$ \cite{dk}.

Suppose we are working over a local ring $(\T, \fm)$.  An important assumption that occurs throughout the literature is that the representations $\rho_i$ are residually distinguishable, i.e.\ that 
\begin{equation} \label{e:rd}
 \overline{\rho}_1 \not\cong \overline{\rho}_2 \pmod{\fm}. \end{equation}
The seminal book on Ribet's method by Bellaiche--Chenevier \cite{bc}, which works in the language of pseudorepresentations, 
assumes throughout that the pseudorepresentations considered are ``residually multiplicity free", i.e.~that the pseudorepresentation analog of (\ref{e:rd}) holds.
The main theorem of this paper is a very general version of Ribet's Lemma for $\GL_2$ where we do not impose the assumption that the associated characters are residually distinguished. The following is a simplified version of our main theorem.

\begin{theorem} \label{t:rl} Let $\T$ be a complete reduced Noetherian local ring,  let $I \subset \T$ be an ideal, and let $\cG$ be a compact group.  Suppose we have a continuous representation
\[ \rho \colon \cG \longrightarrow \GL_2(K), \quad K=\Frac(\T), \]
such that the characteristic polynomial $P_{\rho(g)}(x)$ lies in $\T[x]$ for all $g \in \cG$ and furthermore that
\begin{equation} \label{e:redi}
 P_{\rho(g)}(x) \equiv (x - \chi(g))(x - \psi(g)) \pmod{I} 
 \end{equation}
for two characters $\chi, \psi \colon \cG \longrightarrow \T^*$.
  Suppose  that for every projection of $K$ onto one of its  field factors $K \rightarrow k$, the projection of $\rho$ to 
  $\GL_2(k)$ is an irreducible representation of $\cG$ over $k$.
Then there exists a finitely generated $\T$-module $M$ and a surjective continuous cohomology class \[ \kappa \in H^1(\cG, M(\chi\psi^{-1})) \] such that
\begin{equation} 
 \Fitt_{\T}(M) \subset I.
 \end{equation}
\end{theorem}
Here $M(\chi\psi^{-1})$ denotes the $\T$-module $M$ endowed with a $\cG$-action via $\chi\psi^{-1}$, and
$\Fitt_{\T}(M)$ denotes the 0th Fitting ideal of $M$ over $\T$.
A cohomology class is called surjective if for every representative cocycle $\kappa$, the elements $\kappa(g)$ for $g \in \cG$ generate $M$ as a $\T$-module.   

The novelty of Theorem~\ref{t:rl} is that we do not assume that $\chi \not\equiv \psi \pmod{\fm}$.  
The residually indistinguishable case $\chi \equiv \psi \pmod{\fm}$ was studied recently by Ophir and Weiss, who obtained a version of Ribet's Lemma when working with representations over a DVR \cite{ow}.  As shown by Hajjar Mu\~noz in \cite{rafah}, one can deduce Theorem~\ref{t:rl} when $\T$ is a DVR from the results of Ophir--Weiss.  However, in many arithmetic applications, the local rings $\T$ that occur are Hecke algebras that are rarely DVRs.  In this paper we establish a new technique to handle the  general case. 

\bigskip

Before discussing the proof of Theorem~\ref{t:rl}, we indicate some arithmetic applications of our results. 
As mentioned above, in earlier work the first two named authors proved the Brumer--Stark conjecture away from $p=2$.  Using the main theorem of the present paper, in a companion paper \cite{bsapet} we finish the proof of the Brumer--Stark conjecture by handling the localization at $p=2$.  Let us recall the statement of this conjecture. 

 Let $F$ be a totally real field  of degree $n$ and let $H$ be a finite abelian extension of $F$  that  is a CM field.  Write $G = \Gal(H/F)$.   Let $S$ and $T$ denote finite nonempty disjoint sets of places of $F$ such that $S$ contains the set $S_{\infty}$ of real places and the set $S_{\ram}$ of finite primes ramifying in $H$. 
Associated to any character $\chi \colon G \longrightarrow \C^*$  one has the Artin $L$-function
\begin{equation} \label{e:ls} 
L_S(\chi, s) = \prod_{\fp \not \in S} \frac{1}{1 - \chi(\fp) \N\fp^{-s}}, \qquad \real(s) > 1, 
\end{equation}
and its ``$T$-smoothed" version
\begin{equation} \label{e:lst}
L_{S,T}(\chi, s) = L_S(\chi, s) \prod_{\fp \in T} (1 - \chi(\fp)\N\fp^{1-s}). 
\end{equation}

Assume that $T$ satisfies the Deligne--Ribet condition ensuring the integrality of $L_{S,T}(\chi, 0)$, namely that $T$ contains two  primes of different residue characteristic, or one prime of residue characteristic larger than $n+1$.  In \cite{bsapet} we prove:
\begin{theorem} [Brumer--Stark Conjecture] \label{t:bs}  
 Let $\fp \not\in S \cup T$ be a prime of $F$ that splits completely in $H$.  Fix a prime $\fP$ of $H$ above $\fp$.  There exists an element  $u \in H^*$ 
 satisfying the following.
 \begin{itemize}
\item  We have $|u|_w = 1$ for all places $w$ of $H$ not lying above $\fp$, including the complex places. 
\item We have
\begin{equation} \label{e:bs}
\sum_{\sigma \in G} \chi(\sigma) \ord_{\sigma^{-1}(\fP)}(u) = L_{S,T}(\chi, 0)
 \end{equation}
for all $\chi\colon G \longrightarrow \C^*$. 
\item We have $u \equiv 1 \pmod{\fq\cO_H}$ for all $\fq \in T$.
\end{itemize}
\end{theorem}

In fact, in \cite{bsapet}
we obtain a strong refinement of the Brumer--Stark conjecture that yields the Fitting ideal of certain Ritter--Weiss modules.
In \cite{etnc}, we show how to use this result to deduce
 the minus part of the Equivariant Tamagawa Number Conjecture (ETNC) for the Tate motive associated to $H/F$. Again, we obtain this integrally over $\Z$ and not just over $\Z[1/2]$.  The proof of ETNC, which is obtained by applying an idea of Bullach--Burns--Daoud--Seo \cite{bbds} to our results on Brumer--Stark, yields many important new corollaries, including Rubin's higher rank Brumer--Stark conjecture, the integral Gross--Stark conjecture and higher rank version due to Popescu.  Our results should also yield a version of the classical Main Conjecture of Iwasawa Theory over totally real fields at the prime $p=2$.  Our proof of ETNC is rather formal, with the main arithmetic input arising from the results of~\cite{bsapet}, which in turn are deduced from the main theorem of this paper.

\bigskip
We conclude the introduction by describing some features of the proof of Theorem~\ref{t:rl}.  It is illuminating to 
first consider the residually distinguishable case $\chi \not\equiv \psi \pmod{\fm}$.  In this setting, Theorem~\ref{t:rl} can be proven following Mazur--Wiles (see \cite[Chapter 5, \S 5, Proposition 1]{mw}).  Fix $\tau \in \cG$ such that $\chi(\tau) \not\equiv \psi(\tau) \pmod{\fm}$.  By Hensel's lemma, $P_{\rho(\tau)}(x)$ has two distinct roots in $\T$, 
congruent to $\chi(\tau)$ and $\psi(\tau) \pmod{\fm}$ respectively.
Choose a basis for $\rho$ consisting of the associated eigenvectors for $\rho(\tau)$.  Write \[ \rho(\sigma) = \mat{a(\sigma)}{b(\sigma)}{c(\sigma)}{d(\sigma)} \] and let $B$ denote the $\T$-module generated by the elements $b(\sigma) \in K$ for $\sigma \in G$.  The function 
\[\kappa(\sigma) = \psi^{-1}(\sigma) \overline{b}(\sigma) \in B/IB \]
defines a continuous cocycle representing a surjective class in $H^1(G, B/IB(\chi \psi^{-1}))$.  Furthermore, the irreducibility assumption on $\rho$ implies that $B$ is a faithful $\T$-module which in turn implies that $\Fitt_\T(B/IB) \subset I$.  This concludes our sketch of  the proof in the residually distinguished case.

In the residually indistinguishable case, we define another canonical $\T$-module $M$ and a surjective class $\kappa \in H^1(G, M(\chi\psi^{-1}))$.  The proof that $\Fitt_{\T}(M) \subset I$ is rather elaborate and is the most important contribution of this paper.  New techniques that we introduce to prove this inclusion are the application of {\em matrix invariant theory} and {\em rational cohomology}.  Let us describe these ingredients in greater detail.

By definition, $\Fitt_{\T}(M)$ is the ideal generated by the determinants of all square matrices $D$ of relations occurring in a finite presentation of $M$ over $\T$.
In \S\ref{s:ti} we prove that  assumption (\ref{e:redi}) implies that certain expressions involving the traces and determinants of the matrices $\rho(g)$ lie in the ideal $I$.  It therefore suffices to prove that for all $D$, the element $\det(D) \in \T$ can be expressed in terms of these traces and determinants.  While this is possible to show ``by hand" using some combinatorics in small situations (see for instance \S\ref{s:example}), a general proof requires a more conceptual approach.  

We first pass to a ring of formal variables $R$ that is naturally endowed with a homomorphism $\pi \colon R \lra K$ and prove that it suffices to show that certain expressions in $R$ can be expressed modulo $\ker(\pi)$ in terms of traces and determinants of certain matrices taking values in $R$.  The advantage of this is that the ring $R$ is naturally endowed with an action of the algebraic group  $\G = \GL_2$ over $\Z$.  We then apply a theorem of De Concini  and Procesi, known as the fundamental theorem of matrix invariant theory, that identifies the subring $A \subset R$ generated by traces and determinants as the subring of $R$ invariant under the action of $\G$. 

Now, the element  $x \in R$ naturally mapping to $\det(D) \in \T \subset K$ is not invariant under the action of $\G$, but we prove that its image in $R/J$ is invariant for an appropriate $\G$-invariant ideal $J \subset \ker(\pi)$.  In order to apply the theorem of De Concini and Procesi, we must then show that $x$ lies in the image of the natural map
\[ A = R^{\G} \lra (R/J)^{\G}. \]
The cohomology theory of algebraic group actions goes by the name {\em rational cohomology} (or {\em Hochschild cohomology}).  Here ``rational" refers to actions via rational maps rather than the rational numbers---indeed, it is essential that  we work integrally over $\Z$.  Our goal  is to show that the class $\beta \in H^1(\G, J)$ associated to $x$  under the connecting homomorphism $(R/J)^{\G} \lra H^1(\G, J)$ vanishes.

To do this, we first recall an important result in rational cohomology: the restriction to the lower triangular Borel $\B \subset \G$ induces an isomorphism 
\[ H^i(\G, J) \cong H^i(\B, J) \quad \text{ for all } i \ge 0. \]
Hence it suffices to prove the vanishing of the image of $\beta$ in  $H^1(\B, J)$.  For this, we define a $\B$-invariant subideal $J' \subset J$ such that $\beta$ is in the image of the canonical push-forward  
\begin{equation}
\label{e:pf}
 H^1(\B, J') \lra H^1(\B, J). \end{equation}
We prove that  (\ref{e:pf}) is actually the zero map, and hence that $\beta$  vanishes.  We stress that the ideal $J'$ is only endowed with an action of $\B$ (not the full group $\G$), so restriction to the Borel is a crucial part of our argument.  If the algebraic groups $\B$ and $\G$ are replaced by their $\Z$-valued points, then
the restriction map $H^1(\G(\Z), J) \lra H^1(\B(\Z), J)$ on ordinary group cohomology is in general {\em not} an isomorphism (or even injective).  Our use of rational cohomology in place of group cohomology is therefore essential.

The proof of the vanishing of (\ref{e:pf}) is intricate.  We begin by constructing a resolution $C^\bullet$ of the ideal $J'$ by $\B$-modules. Here
  the Kozsul complex giving a resolution of an ideal generated by a regular sequence and its generalization by Buchsbaum--Rim to determinantal ideals play an important role.  
  Next we embed $C^\bullet$ into a complex $D^\bullet$ such that $H^0(D^\bullet) = R/J$, and such that $D^\bullet$ consists of acyclic $\B$-modules.
We then prove a general result showing that in such a setup (i.e.\ where $C^\bullet$ is bounded and exact, and $D^\bullet$ is acyclic), the push forward (\ref{e:pf}) necessarily vanishes. See Theorems \ref{t:push} and \ref{t:comm} for precise statements.  
  
  \bigskip
We conclude the introduction by noting that in the main text we actually prove a stronger version of Theorem~\ref{t:rl} that  relaxes the assumption that $\T$ is reduced and that establishes local conditions for the cohomology class we construct.  See Theorem~\ref{t:rll} below.
It is the local conditions that require the use of the Buchsbaum--Rim complex; for the global picture in Theorem~\ref{t:rl}, the classical Koszul complex suffices.  The reader who wishes to understand the main aspects of our argument in a simplified setting (e.g.\ where we ignore local conditions) is encouraged to consult the announcement \cite{das}.

It is natural to attempt to generalize our construction beyond the setting of $\GL_2$.  Our hope is that our construction in the residually indistinguishable case of Ribet's method, appropriately generalized, will have arithmetic applications beyond those presented in \cite{bsapet} and described above. 

\bigskip

We would like to thank Rafah  Hajjar Mu\~noz, who visited the first named author at Duke University from UPC (Barcelona) to write an undergraduate senior thesis in the 2021--2022 academic year.  Rafah did computer examples generalizing \S\ref{s:example}, and it was by analyzing the formulas he produced that we realized the role of the Koszul complex in the proof of  the vanishing of $\beta \in H^1(\G, J)$.  We would also like to thank Robert Boltje, Brian Conrad, Corrado De Concini, Claudio Procesi, and Geordie Williamson for helpful discussions.

The first named author is supported by a grant from the National Science Foundation (DMS-2200787). The second named author is supported by DST-SERB grant SB/SJF/2020-21/11, SERB MATRICS grant MTR/2020/000215, SERB SUPRA grant SPR/2019/ 000422, and DST FIST program - 2021 [TPN - 700661]. The fourth named author is supported by NSF grant DMS 2201346.

\section{Main Theorem and Construction of Cocycle}\label{s:cocycle}

Our main theorem is a strengthening of Theorem~\ref{t:rl} stated in the introduction, in which we  incorporate certain local conditions at the places of $F$.  Local conditions are always necessary in arithmetic applications, and we have chosen conditions here tailored to the application for Brumer--Stark in  \cite{bsapet}.  We also relax the assumption that $\T$ is reduced, though we still need an assumption on the ring $K$ in which our representation $\rho$ lands.

\begin{theorem} \label{t:rll}
 Let $\T \subset \tilde{\T}$ be an inclusion of commutative Noetherian rings, with $\T$ local.  Suppose that $\T$ and $\tilde{\T}$ are complete with respect to the maximal ideal of $\T$.  Let $\tilde{I} \subset \tilde{\T}$ be a nontrivial ideal and let $I = \tilde{I} \cap \T$.  Let $K = \Frac(\tilde{\T})$ be the total ring of fractions of $\tilde{\T}$.
 Assume that $K$ is a product of local rings and that the maximal ideals of $K$ are principal.
  Let $\cG$ be a compact group.
 Suppose we are given a continuous  representation
\[ \rho \colon \cG \longrightarrow \GL_2(K) 
\]
satisfying the following conditions.
\begin{itemize}
\item  
 For each $\sigma \in \cG,$  the characteristic polynomial $P_{\rho(\sigma)}(x)$ lies in $\T[x]$.  Furthermore we have
\begin{equation} \label{e:redpoly}
P_{\rho(\sigma)}(x) \equiv (x - \chi(\sigma))(x - \psi(\sigma)) \pmod{I} 
\end{equation}
for two characters $\chi, \psi \colon \cG \longrightarrow \T^*$ such that $\chi \equiv \psi \pmod{\fm}$.
\item   Let $K_0 = 
\red(K)$ denote the maximal reduced quotient of $K$. Write $K_0 = \prod_{i=1}^{m} k_i$ as a product of fields.  For every projection $K \rightarrow K_0 \rightarrow k_i$, the projection of $\rho$ to $\GL_2(k_i)$ is an irreducible representation of $\cG$ over $k_i$.
\item We are given a  set of subgroups $\cG_v \subset \cG$, indexed by a finite set $S$, such that for each $v \in S$ there exists a basis in which the restriction of $\rho$ has the form
\begin{equation} \label{e:rholocal}
 \rho|_{\cG_v} \cong \mat{\eta_v}{0}{*}{\xi_v}
\end{equation}
for two characters $\xi_v, \eta_v \colon \cG_v \longrightarrow \tilde{\T}^*.$  
\item We are given a partition $S = \Sigma \sqcup \cP$.  For each $v \in \Sigma$ we have the congruence 
  \begin{equation} \xi_v \equiv \psi|_{\cG_v} \pmod{\tilde{I}}. \label{e:Sigmaassumption}
   \end{equation}
\item For each $v \in \cP$ we are given a subgroup $\cI_v \subset \cG_v$ such that 
\begin{equation} \label{e:Passumption}
(\xi_v)|_{\cI_v} \equiv \chi|_{\cI_v} \pmod{\tilde{I}}. \end{equation}
   \end{itemize}
 If $\Sigma$ is nonempty,  fix $v_0 \in \Sigma$.  Choose an element $\sigma_v \in \cG_v$ for each $v \in \cP$.  Then there exists a finitely generated $\T$-module $N$ and a 
continuous cocycle \[ \kappa \in Z^1(\cG, N(\chi\psi^{-1})) \] 
satisfying the following conditions.
\begin{itemize}
\item If $\Sigma$ is nonempty, we have $\kappa(\cG_{v_0}) = 0$ and for each $v \in \Sigma \setminus \{v_0\}$, there exists $y_v \in N$ such that
\[ \kappa(\sigma) = (\chi\psi^{-1}(\sigma) - 1) y_v \]
for all $\sigma \in \cG_v.$
\item The module $N(\chi\psi^{-1})$ is generated over $\T$ by $\kappa(\cG)$ and the $y_v$, $v \in \Sigma \setminus \{v_0\}$.
\item For each $v \in \cP$, we have $\kappa(\sigma) = 0$ for all $\sigma \in \cI_v$.

\item We have \begin{equation} \label{e:nincl}
 \left(\prod_{v \in \cP} (\xi_v(\sigma_v) - \chi(\sigma_v))\right)\Fitt_{\T}(N) \subset \tilde{I}. \end{equation}
\end{itemize}
\end{theorem}

\begin{remark} We have included the assumption of residual indistinguishability \[ \chi \equiv \psi \pmod{\fm} \] in the statement of Theorem~\ref{t:rll} since that is the salient case for this paper.  The theorem  remains true without that assumption and can be proven when $\chi \not\equiv \psi \pmod{\fm}$ using the ``$b$-coefficient'' of the representation $\rho$ in the appropriate basis, as indicated in the introduction.  The theorem was essentially proven this way in the residually distinguishable case in \cite{dk}, though it was not stated in  this precise form.
\end{remark}

\subsection{Construction of the module $N$}

Our setting is as in the statement of Theorem~\ref{t:rll}.
Extend $\rho$ to a continuous $\T$-algebra homomorphism
\[ 
\T[\cG] \longrightarrow M_2(K). 
\] 
Similarly extend $\chi, \psi$ to $\T$-algebra homomorphisms $\T[\cG] \longrightarrow \T$.
It is well-known that the congruence (\ref{e:redpoly}) extends to all $t \in \T[\cG]$:  
\begin{equation} \label{e:redpoly2}
P_{\rho(t)}(x) \equiv (x - \chi(t))(x - \psi(t)) \pmod{I}.
\end{equation}
See for instance \cite{ow}*{Lemma 3.1}.

Define two $\T$-submodules of $M_2(K)$:
 \[ \Delta_\chi = \{ \rho(t) - \chi(t) \colon t \in \T[\cG]\}, \qquad \Delta_{\psi} =   \{ \rho( t ) - \psi(t)\colon t \in \T[\cG]\}. \]
  Here $\chi(t)$ and $\psi(t)$ denote scalar matrices.
It is elementary to check that $\Delta_{\chi} \Delta_{\psi} \subset \Delta_{\psi}$, where the module on the left represents the $\T$-module generated by all products $tt'$ with $t \in \Delta_\chi, t' \in \Delta_\psi$.    We then define the $\T$-module
\begin{equation} M_0 = \Delta_{\psi}/ \Delta_{\chi} \Delta_{\psi}.
\end{equation}

\begin{lemma} The map $\kappa \colon \cG \longrightarrow \Delta_{\psi}/\Delta_{\chi}\Delta_{\psi}$ given by
\[
 g \mapsto \psi^{-1}(g)(\rho(g) - \psi(g))
\]
defines a continuous cocycle $\kappa \in Z^1(\cG, M_0(\chi\psi^{-1}))$.
\end{lemma}
\begin{proof}
We check the cocycle condition:
\begin{align*}
\kappa(g_1g_2) \ - & \ \kappa(g_1) - \psi^{-1}\chi(g_1)\kappa(g_2) \\
= & \ \kappa(g_1g_2) - \psi^{-1}(g_1)(\rho(g_1) - \psi(g_1)) - \psi^{-1}\chi(g_1) \psi^{-1}(g_2)(\rho(g_2) - \psi(g_2)) \\
= & \  \psi^{-1}(g_1g_2)(\rho(g_1g_2) - \psi(g_1g_2)) \\
& \ \ \ \ \  - \psi^{-1}(g_1) \rho(g_1) - \psi^{-1}(g_1g_2)\chi(g_1) \rho(g_2) - \psi^{-1}(g_1)\chi(g_1) \psi(g_2) + 1 \\ 
= & \  \psi^{-1}(g_1g_2)(\rho(g_1) - \chi(g_1))(\rho(g_2) - \psi(g_2)).
\end{align*}
The last item vanishes in the quotient defining $M_0$.
\end{proof}  

Define \[ N_0 = M_0 \oplus \bigoplus_{v \in \Sigma \setminus \{v_0\}} \T y_v, \]
i.e.\ the direct sum of $M_0$ with the free $\T$-module on the set $\Sigma \setminus \{v_0\}$. 
Let $Q \subset N_0$ denote the $\T$-submodule generated by the relations we must impose for $\kappa$ to satisfy the desired properties, namely:
\begin{itemize}
\item $\kappa(\sigma)$ for $\sigma \in \cG_{v_0}$ (in the case that $\Sigma$ is nonempty).
\item $\kappa(\sigma) - (1 - \chi \psi^{-1}(\sigma)) y_v $ for $\sigma \in \cG_v$, $v \in \Sigma \setminus \{v_0\}$. 
\item $\kappa(\sigma)$ for $\sigma \in \cI_v$, $v \in \cP$.
\end{itemize}
We then define
\[ N = N_0/Q \]
and let $M$ be the image of $M_0$ in $N$.
Note that $N$ is finitely generated since $\cG$ is compact, $\rho$ is continuous, and $\Sigma$ is finite.  The first three bullet points required of the module $N$ in Theorem~\ref{t:rll} are clearly satisfied.

The remainder of the paper is taken up in proving the last bullet point of Theorem~\ref{t:rll}, which states
\begin{equation} \label{t:lbp}
 \left(\prod_{v \in \cP} (\xi_v(\sigma_v) - \chi(\sigma_v))\right)\Fitt_{\T}(N) \subset \tilde{I}.  
 \end{equation}
  
\subsection{Fitting Ideal} \label{s:fi}

If $\Sigma$ is nonempty, we fix the basis for the representation $\rho$ corresponding to (\ref{e:rholocal}) for the place $v_0$, i.e.\ we choose a basis such that $\rho|_{\cG_{v_0}}$ is lower triangular. If $\Sigma$ is empty, any basis will suffice.
Write \[ \rho(\sigma)- \psi(\sigma) = \mat{a(\sigma)}{b(\sigma)}{c(\sigma)}{d(\sigma)}. \]
Let \begin{equation} \label{e:rhoidef}
 \left\{\rho_i = \rho(g_i) - \psi(g_i) = \mat{a_i }{b_i}{c_i}{d_i}: 1 \leq i \leq r\right\},
\end{equation}
 with $g_i \in \cG$, be a set of $\T$-module generators for $\Delta_{\psi}$ such that $ \left\{ \rho(g_i) - \chi(g_i) \right \} $ is a set of $\T$-module generators for $\Delta_\chi$.
The module $N$ is generated over $\T$ by the images of the $\rho_i$ together with the images of the $y_v$ for $v \in \Sigma \setminus \{v_0\}$.  Write $s = \#\Sigma \setminus \{v_0\}$ (and $s=0$ if $\Sigma$ is empty).
For each $\sigma \in \cG$ choose coefficients $\alpha_{\sigma, i} \in \T$ such that
\[ \rho(\sigma) - \psi(\sigma) = \sum_{i=1}^r \alpha_{\sigma,i} \rho_i. \]

 There are 5 types of relations among the $r+s$ generators $\{\rho_i \} \cup \{ y_v \}$ in the module $N$:
\begin{itemize}
\item[(I)] Relations among the $\rho_i$ that already hold in $\Delta_\psi$, before any quotient is taken.  Suppose
\begin{equation} \label{e:epsilondef}
\sum_{k=1}^r \epsilon_k \rho_k = 0  \qquad \text{ with } \epsilon_k \in \T.
\end{equation}
We have the relation
\begin{equation} \label{e:epsrel}
(\epsilon_1, \dots, \epsilon_r, 0, \dots, 0)  \qquad (s \text{ zeros}).
\end{equation}
\item[(II)]  Relations arising from the quotient by $\Delta_\chi \Delta_\psi$.  Write 
\begin{equation} \label{e:nudef}
 \nu_i = \psi(g_i) - \chi(g_i) \end{equation} 
with the $g_i$ as in (\ref{e:rhoidef}). 
  For each $1 \le i, j \le r$, write
\begin{equation} \label{e:deltadef}
(\rho_i + \nu_i) \rho_j = \sum_{k = 1}^r \delta_{ijk} \rho_k
\end{equation}
 with $\delta_{ijk} \in \T$. 
We have the relation
\begin{equation} \label{e:delrel}
(\delta_{ij1}, \dots, \delta_{ijr}, 0, \dots, 0). 
\end{equation}
\item[(III)] Relations arising from the quotient by $\kappa(\cG_{v_0})$.  For each $\sigma \in \cG_{v_0}$, we have the relation
\begin{equation} \label{e:alpharel}
(\alpha_{\sigma,1}, \dots, \alpha_{\sigma,r}, 0, \dots, 0). 
\end{equation}
\item[(IV)] Relations arising from the quotient by $\kappa(\cI_v)$ for $v \in \cP$.  For each $\sigma \in \cI_v$, we have the relation
\begin{equation} \label{e:betarel}
(\alpha_{\sigma,1}, \dots, \alpha_{\sigma,r}, 0, \dots, 0). 
\end{equation}
\item[(V)] Relations arising from the quotient by the conditions at $v \in \Sigma$.  For each $\sigma \in \cG_v$, we have the relation
\begin{equation} \label{e:gammarel}
(\alpha_{\sigma,1}, \dots, \alpha_{\sigma,r}, 0, \dots, 0, \psi(\sigma) - \chi(\sigma), 0). 
\end{equation}
Here the entry $\psi(\sigma) - \chi(\sigma)$ occurs in the coordinate corresponding to $y_v$.
\end{itemize}

By definition, $\Fitt_{\T}(N)$ is the ideal of $\T$ generated by the determinants of all square matrices $D_0$ of dimension $(r+s)$ whose rows are  vectors of any of the forms (\ref{e:epsrel}), (\ref{e:delrel}), (\ref{e:alpharel}), (\ref{e:betarel}), or (\ref{e:gammarel}). 

\subsection{Notational simplification} \label{s:ns}

Suppose we are given a matrix $D_0$ as above.
In what follows, it is notationally convenient if we replace our given matrix $D_0$ with another matrix $D$ that has the same determinant.
Namely, we may choose a larger set of generators of $\Delta_{\psi}$ so that, for all $\sigma \in  \cG$ contributing rows of types (III), (IV), or (V) to $D_0$, we include as one of our generators $\rho_i$ the element 
$\rho_\sigma = \rho(\sigma) - \psi(\sigma)$.

This change replaces the matrix $D_0$ with a matrix $D$ defined as follows:
\begin{itemize} 
\item The new matrix $D$ has one additional column associated to each row of type (III), (IV), or (V) in $D_0$, i.e. to each new generator $\rho_\sigma$.
\item Each row of type (I) or (II) in $D_0$ has a corresponding row in $D$, with 0's in the new columns.
\item  Each row of type (III), (IV), or (V) in $D_0$ gets replaced by 2 rows in $D$:
\begin{itemize}
\item One row $(\alpha_{\sigma,1}, \dots, \alpha_{\sigma,r}, 0, \dots, 0, -1, 0 , 0 , \dots, 0)$, with the $-1$ in the new column corresponding to $\rho_\sigma$. Note that this is a  row of type (I).
\item One row \begin{equation} \label{e:newtype}
(0, 0 , \dots, 1, 0 , \dots , 0), \end{equation}
 with the $1$ in the $\rho_\sigma$ column if we have a row of type (III) or (IV), or 
one row \begin{equation} \label{e:newtypeV} (0,  0 , \dots, 1, 0 , \dots, 0, \psi(\sigma) - \chi(\sigma), 0, \dots, 0), \end{equation} with the $1$ in the $\rho_\sigma$ column and $\psi(\sigma) - \chi(\sigma)$ in the $y_v$ column if we have a row of type (V).
\end{itemize}
\end{itemize}

 It is easy to see that the new matrix $D$ has the same determinant as the old one, but the form of the rows of type (III), (IV), or (V) has been greatly simplified.
Going forward, we will ignore $D_0$ and only work with the matrix $D$. In particular, $r$ denotes the number of generators corresponding to $D$ (i.e., the $r$ for $D_0$ plus the number of new generators included).  For the matrix $D$, the rows of type (III) or (IV) all have the form (\ref{e:newtype}), and the rows of type (V) all have the form (\ref{e:newtypeV}).

\subsection{Auxiliary Matrices}  Write $\cP = \{ v_1, \dots, v_t\}.$ For notational simplicity, write $\sigma_i = \sigma_{v_i}$ and define
\begin{equation} \label{e:zidef}
z_{\sigma_{i}} = \xi_{v_i}(\sigma_{i}) - \chi(\sigma_{i}). 
\end{equation}
  Associated to our matrix $D$ from \S\ref{s:ns}, we define an auxiliary $(r + s + t) \times (r+s + t)$ matrix $E$ with coefficients in $\T$ as follows. 

\begin{equation} E =  \left( \ 
\begin{tikzcd}[cramped, sep=small, row sep=-3pt]
z_{\sigma_1} &0 &[-30pt] \cdots  &[-30pt]  0 & \alpha_{\sigma_{1},1} & \cdots & \alpha_{\sigma_{1},r} & 0 & \cdots & 0 \\
0  & z_{\sigma_2} & \cdots & 0 & \alpha_{\sigma_{2},1} & \cdots & \alpha_{\sigma_{2},r} & 0 & \cdots & 0 \\
    &      & \ddots &  &   & \vdots &   &  & \vdots &  \\
    0  & 0 & \cdots &  z_{\sigma_t} & \alpha_{\sigma_t,1} & \cdots & \alpha_{\sigma_t,r} & 0 & \cdots & 0 \\
 \ \\
 \ \\
    &     &  \!\!\!\!\! \scaleto{0_{(r+s)\times t}}{24pt} & & & &  \scaleto{D}{24pt} & & & 
 \end{tikzcd} \  \right)
 \end{equation}

\bigskip

Clearly we have $\det(E) = (\prod_{i =1}^t z_{\sigma_i}) \det(D)$, so our goal in proving (\ref{t:lbp}) is to show that $\det(E) \in \tilde{I}$. 
 To prove this, we will define an alteration of the matrix $E$ yielding a new square matrix $E'$ of dimension $(r+s+t)$ with coefficients in $K = \Frac(\tilde{\T})$, rather than $\tilde{\T}$.  The motivation behind these alterations will become clear in Lemma~\ref{l:detzero} below.
For the first $t$ rows of $E$, we define
\begin{equation} \label{e:xidef}
x_{\sigma_{i}} = \xi_{v_i}(\sigma_{i}) - \psi(\sigma_{i})
\end{equation}
and we  replace the element $z_{\sigma_i}$ by 
\[ x_{\sigma_i} - a(\sigma_{i}) = x_{\sigma_i} - \sum_{j=1}^{r} a_j \alpha_{\sigma_i,j}. \]

For the last $r+s$ rows of $E$, we alter the rows as follows, based on the type of relation appearing in the corresponding row of $D$:
\begin{itemize}
\item[(I)] or (III) No change.
\item[(II)] Replace $\delta_{ijj}$ by $\delta_{ijj} - a_i - \nu_i$, replace $\delta_{iji}$ by $\delta_{iji} - d_j$, and leave the other $\delta_{ijk}$ unchanged.
\item[(IV)] Insert \[ \xi_v(\sigma) - \psi(\sigma) - a(\sigma) \] into the column among the first $t$ columns of $E$ corresponding to $v$.
\item[(V)] Replace $\psi(\sigma) - \chi(\sigma)$ by 
\[ \xi_v(\sigma) - \psi(\sigma) - a(\sigma). \]
\end{itemize}

The purpose of these replacements is that with these small changes, the matrix $E'$ obtained has vanishing determinant.

\begin{lemma} We have $\det(E') = 0.$ \label{l:detzero}
\end{lemma}

\begin{proof}  By assumption, the ring $K$ is a product of local rings, and the maximal ideal of each of these local rings is principal.  It suffices to prove the result on each such local factor, so we assume that $K$ is a local ring with maximal ideal generated by an element $\pi$.  We will prove the result by demonstrating an appropriate vector in the kernel of $E'$.

We are given that for each place $v \in S$, there exists a change of basis matrix \[ M_v = \mat{A_v}{B_v}{C_v}{D_v} \in \GL_2(K) \] such that
\begin{equation} \label{e:changebasis}
  \mat{a(\sigma)+ \psi(\sigma)}{b(\sigma)}{c(\sigma)}{d(\sigma) + \psi(\sigma)}M_v = M_v\mat{\eta_v(\sigma)}{0}{*}{\xi_v(\sigma)}  \end{equation}
for all $\sigma \in \cG_v$.  Equating the upper right entries of this matrix equation, we obtain 
\begin{equation} \label{e:db}
 D_v b(\sigma) = B_v(\xi_v(\sigma) - \psi(\sigma) - a(\sigma)). 
 \end{equation} 
Let us for the moment assume that each $D_v$ is invertible in $K$.  
Writing $\cP = \{v_1, \dots, v_t\} $ and $\Sigma \setminus \{v_0\} = \{w_1, \dots, w_s \}$, 
we define a (column) vector of length $r+s+t$ as follows:
\[ w = \left(- \frac{B_{v_1}}{D_{v_1}}, \dotsc, -\frac{B_{v_t}}{D_{v_t}}, b_1, b_2, \dots, b_r, -\frac{B_{w_1}}{D_{w_1}}, \dotsc, -\frac{B_{w_s}}{D_{w_s}}\right). \]
The fact that $E' w = 0$ follows from the definition of the alterations made in the definition of $E'$.  For example, for the first $t$ rows, the dot product of the $i$th row of $E'$ with the vector $w$ is
\begin{align*}
& \ (\xi_{v_i}(\sigma_{i}) - \psi(\sigma_{i}) - a(\sigma_{i}))\left(-\frac{B_{v_i}}{D_{v_i}}\right) + \sum_{j=1}^r \alpha_{\sigma_i,j} b_j  \\
   = & \  (\xi_{v_i}(\sigma_{i}) - \psi(\sigma_{i}) - a(\sigma_{i}))\left(-\frac{B_{v_i}}{D_{v_i}}\right) + b(\sigma_{i}) \\
  = & \ 0 
  \end{align*}
    by (\ref{e:db}). 
 The orthogonality with the other rows is similar, and we content ourselves with brief explanations: 
  \begin{itemize}
  \item For rows of type IV or V, we again use (\ref{e:db}).
  \item For rows of type I, we use
  \[ \sum_{j=1}^r \epsilon_{ij} b_j = 0. \]
The same  holds for rows of type (III)---this is where our chosen basis for $\rho$ is used (recall we chose the basis such that $\rho|_{\cG_{v_0}} $ is lower triangular).
\item For rows of type II, we use the equation
\[ (a_i+\nu_i)b_j + b_i d_j = \sum_{k=1}^r \delta_{ijk} b_k \]
arising from the ``$b$''-component of the equation (\ref{e:deltadef}) defining the $\delta_{ijk}$.
  \end{itemize}
  The second bulleted assumption in the statement of Theorem~\ref{t:rll} regarding the irreducibility of the projections of $\rho$ implies that the elements $b_1, \dots, b_r$ generate $K$ as a $K$-module.  Therefore the equation $E'w = 0$ implies that $\det(E') = 0$ as desired.
  
  This argument is easily adapted in the case that some $D_{v_i}$ is not a unit in $K$.  Recall we have reduced to the setting that $K$ is a local ring with maximal ideal generated by an element $\pi$ (which is necessarily nilpotent).  For each $v \in S$, we can write $D_{v} = u_v \pi^{e_v}$ for some unit $u_v \in K^*$ and integer $e_v \ge 0$.   Let $m$ denote the maximal $e_v$ (if some $D_{v} =0$, then we let 
$m = e_v$ be the minimal exponent $k$ such that $\pi^k = 0$).  Equation (\ref{e:db}) yields
\[  \pi^m b(\sigma) = \frac{B_{v_i}}{u_{v_i}} \pi^{m - e_{v_i}}(\xi_{v_i}(\sigma) - \psi(\sigma) - a(\sigma)). 
\]
We therefore define
\[ w = \left(- \frac{B_{v_1}}{u_{v_1}} \pi^{m-e_{v_1}}, \dotsc, -\frac{B_{v_t}}{u_{v_t}}\pi^{m-e_{v_t}}, b_1 \pi^m, b_2 \pi^m, \dots, b_r \pi^m, -\frac{B_{w_1}}{u_{w_1}}\pi^{m-e_{w_1}}, \dotsc, -\frac{B_{w_s}}{u_{w_s}}\pi^{m-e_{w_t}}\right). \]
The equation $E'w = 0$ follows as before.  Now for the $v$ yielding the maximal $e_v = m$, the component of $w$ is 
$-B_v/u_v \in K^*$.  Here we use the fact that $D_v \not\in K^* \Longrightarrow B_v \in K^*$ since $M_v \in \GL_2(K)$.  Therefore $E' w = 0$ again implies $\det(E') = 0$.
\end{proof}

Our goal is to show that $\det(E) \in \tilde{I}$, and we have shown that $\det(E') = 0$.  It therefore suffices to show that 
$\det(E') - \det(E) \in \tilde{I}$.  In other words, the alterations used to pass from $E$ to $E'$ are small enough to leave the determinant unchanged modulo $\tilde{I}$. Let us motivate our strategy to prove this with an example.

\subsection{An example} \label{s:example}

Suppose $r=2$ and that $S = \Sigma \sqcup \cP$ is empty.   We consider a matrix with only rows of type II, namely

\[ E = D = \mat{\delta_{121}}{\delta_{122}}{\delta_{211}}{\delta_{212}}, \qquad \text{ whence } \qquad E' = \mat{\delta_{121} - d_2}{\delta_{122} - a_1 - \nu_1}{\delta_{211} - a_2 - \nu_2}{\delta_{212} - d_1}. \]
As shorthand, write $t_i = a_i + d_i$ for the trace of $\rho_i$ and $t_{12}$ for the trace of $\rho_1 \rho_2$.
By multilinearity of the determinant, we have \begin{align}
\det(E') - \det(E) &= - \det\mat{d_2}{a_1+\nu_1}{\delta_{211}}{\delta_{212}} - \det \mat{\delta_{121}}{\delta_{122}}{a_2+\nu_2}{d_1}
+ \det \mat{d_2}{a_1 + \nu_1}{a_2+\nu_2}{d_1} \nonumber \\
& = ( t_1 + \nu_1) \delta_{211} - (d_1 \delta_{211} + d_2 \delta_{212})  +   \nonumber \\
& \ \  \ \  + (t_2 + \nu_2) \delta_{122} - (d_1 \delta_{121} + d_2 \delta_{122} )  +   \nonumber\\
&  \ \ \ \   +  (d_1d_2 - (a_1 + \nu_1)(a_2 + \nu_2))   \nonumber \\
\begin{split}
&= ( t_1 + \nu_1) \delta_{211} - (c_2b_1 + (d_2 + \nu_2)d_1)  +   \label{e:exdel}  \\
& \ \  \ \  + (t_2 + \nu_2) \delta_{122} - (c_1b_2 + (d_1 + \nu_1)d_2 )  +   \\
&  \ \ \ \   +  (d_1d_2 - (a_1 + \nu_1)(a_2 + \nu_2))    
\end{split} \\
&= ( t_1 + \nu_1) \delta_{211} + (t_2 + \nu_2) \delta_{122} - (t_{12} + t_1 \nu_2 + t_2 \nu_1 + \nu_1\nu_2 ).
\label{e:example}
 \end{align}
 Equation (\ref{e:exdel}) uses the definition (\ref{e:deltadef}) of the $\delta$'s.
By the fundamental assumption concerning the representation $\rho$ (namely, the first bullet point in the statement of Theorem~\ref{t:rll}), we know that $\tr(\rho(\sigma)) - \psi(\sigma) - \chi(\sigma) \in I$. Recalling the definition of $\nu_i$ from (\ref{e:nudef}), it directly follows that $t_i + \nu_i \in I$.  A similar elementary computation that is generalized below also shows that $t_{12} + t_1 \nu_2 + t_2 \nu_1 + \nu_1\nu_2 \in I$. Hence (\ref{e:example}) lies in $I$ as desired.

\subsection{Trace identities} \label{s:ti}

We now return to the general setting.  We prove that certain expressions involving traces and determinants, generalizing the elements in parentheses in (\ref{e:example}), lie in $I$.

\begin{lemma}\label{l:dets} We have $\det(\Delta_\psi) \subset I$.
\end{lemma}
\begin{proof}
We have $\det(\rho(t) - \psi(t)) = P_{\rho(t)}(\psi(t)) \equiv 0 \pmod{I}$ by (\ref{e:redpoly2}).
\end{proof}

Let $\T\langle X_1, \dots, X_r \rangle$ denote the polynomial algebra over $\T$ in $r$ noncommuting variables $X_i$.  Define a $\T$-algebra homomorphism \[ V \colon \T\langle X_1, \dots, X_r \rangle \longrightarrow \T, \qquad X_i \mapsto -\nu_i. \]

\begin{lemma}\label{l:tr-char}
For any $f \in  \T\langle X_1, \dots, X_r \rangle$ with constant term 0, we have \[ \tr(f(\rho_1, \dotsc, \rho_r)) \equiv V(f) \pmod{I}. \]
\end{lemma}
\begin{proof}
We demonstrate the result for a monomial $X_{j_1}X_{j_2} \cdots X_{j_k}$ with $k \ge 1$.
Write 
\[ \rho_{j_i} = \rho(g_i) - \psi(g_i), \qquad \chi_i = \chi(g_i), \psi_i = \psi(g_i). \]
We calculate the following, where the congruence holds modulo $I$:
\begin{align}
\tr\Big(\prod_i \rho_{j_i}\Big) &= \tr\Big(\prod_i (\rho(g_i) - \psi(g_i))\Big) \nonumber \\
&= \tr\Big(\prod_i \rho(g_i)\Big) - \sum_{i \neq j} \psi_i \tr\Big(\prod_{k \neq i} \rho(g_k)\Big)  \nonumber \\
& \ \ \ \ \  +  \sum_{i \neq j} \psi_i \psi_j \tr\Big(\prod_{k \neq i, j} \rho(g_k)\Big) + \cdots + (-1)^n 2 \prod_i \psi_i  \nonumber \\
\begin{split}
&\equiv \Big(\prod_i \chi_i + \prod_i \psi_i\Big) - \sum_{i} \psi_i \Big(\prod_{k \neq i} \chi_i + \prod_{k \neq i} \psi_i\Big)   \\
 & \ \ \ \ \ + \sum_{i \neq j} \psi_i \psi_j \Big(\prod_{k \neq i, j} \chi_i + \prod_{k \neq i, j} \psi_j\Big) + \cdots + (-1)^n 2\prod_i \psi_i. \label{e:ti}
\end{split}
\end{align}
The $k$-th term has ${n \choose k}$ copies of $\prod_i \psi_i$, except the $(n+1)$st (final) term which has $2$ copies. Taking into account the signs, we see that the total contribution of $\prod_i \psi_i$ to the sum is $(1-1)^n + (-1)^n = (-1)^n$. The expression (\ref{e:ti}) becomes
\begin{align*}
 \prod_i \chi_i & - \sum_{i} \psi_i \prod_{k \neq i} \chi_i  + \sum_{i \neq j} \psi_i \psi_j \prod_{k \neq i, j} \chi_i + \cdots + (-1)^n \prod_i \psi_i \\
 & = \prod_i (\chi_i - \psi_i) = \prod_i (-\nu_{j_i}). 
\end{align*}
This is the desired result.
\end{proof}

\section{Formal Variables} \label{s:formal}

The goal of the remainder of the paper is to prove that
 \[ \det(E') - \det(E) \in \tilde{I} \] as in the example of \S\ref{s:example}.  For this, we must prove the existence of certain algebraic identities relating 
 $\det(E') - \det(E) $ to expressions involving traces and determinants, allowing for the application of the lemmas in \S\ref{s:ti}.  It seems intractable to derive explicit formulas such as (\ref{e:example}) in the general setting.  Instead, a more abstract proof is required.

Since we are aiming to prove the existence of a certain polynomial identity, we shift our perspective from working with the ring $\T$ to working with formal polynomial rings.  We will define a polynomial algebra $R$ and a specialization homomorphism $\pi\colon R \longrightarrow K = \Frac(\tilde{\T})$.  

\subsection{The ring $R$} \label{s:trr}

We first define the ring
\begin{equation} \label{e:r0def}
 R_0 = \Z[ \boldsymbol{\nu}_i, \boldsymbol{\epsilon}_{i}, \boldsymbol{\delta}_{ijk}, \boldsymbol{x}_\sigma   ], 
 \end{equation}
whose generators, along with their images under $\pi$, we specify precisely as follows:
\begin{itemize}
\item Variables $\boldsymbol{\nu}_i$ for $i = 1,\dotsc, r$ with $\pi(\boldsymbol\nu_i) = \nu_i $ as in (\ref{e:nudef}).
\item Variables $\boldsymbol{\epsilon}_{i}, \boldsymbol{\delta}_{ijk}$ corresponding to the rows of type (I) or (II) in the matrix $D$, with $\pi(\boldsymbol{\epsilon}_{i}) = \epsilon_i$ 
and $\pi(\boldsymbol{\delta}_{ijk}) = \delta_{ijk}.$  To be clear, there is one set of variables $\{\boldsymbol{\epsilon}_{1}, \dots, \boldsymbol{\epsilon}_{r}\}$ for each row of type I of $D$, and one set of variables
$\{\boldsymbol{\delta}_{ij1}, \dots, \boldsymbol{\delta}_{ijr}\}$ for each row of type II of $D$ (in the latter case, with $(i,j)$ determined by the row).
\item Variables $\boldsymbol{x}_\sigma$ for each row of type IV or V in $D$ associated to an element $\sigma \in \cG_v$, or for $\sigma = \sigma_i$ from the first $t$ rows of $E$.  We have $\pi(\boldsymbol{x}_\sigma) = \xi_v(\sigma) - \psi(\sigma)$.
\end{itemize}

We then define \begin{equation} \label{e:r1def}
 R_1 = R_0[\boldsymbol{a}_i, \boldsymbol{b}_i, \boldsymbol{c}_i, \boldsymbol{d}_i]_{i=1}^r, 
 \end{equation}
with $\pi(\boldsymbol{a}_i) = a_i$, etc., the components of the matrices $\rho_i$ (see (\ref{e:rhoidef})).  

Finally, recall that $\cB_{v_0}$ denotes the set of $\sigma$ indexing the rows of $D$  of type III, corresponding to certain generators $\rho_{\sigma}$ ($= \rho_i$ for some $i$). 
We then write $\boldsymbol{b}_\sigma = \boldsymbol{b}_i$ and define
\begin{equation} R = R_1/(\boldsymbol{b}_{\sigma} \colon \sigma \in \cB_{v_0}). \label{e:rdef} 
\end{equation}
As $\pi(\boldsymbol{b}_{\sigma}) = 0$, we obtain an induced map $\pi \colon R \lra K.$

\bigskip

We next define matrices $\boldsymbol{E}, \boldsymbol{E}' \in M_{r+s+t}(R)$ whose images under $\pi$ are related to $E, E'$, respectively.  The matrix $\boldsymbol{E}$ is defined from $E$ simply by making every entry bold, with three exceptions. 
\begin{itemize}
\item In the first $t$ rows of $E$ we replace
$z_{\sigma_i}$ by $\boldsymbol{x}_{\sigma_i} + \boldsymbol{\nu_{\sigma_i}}$. 
\item For the rows of type IV, we insert  $\boldsymbol{x}_\sigma + \boldsymbol{\nu_{\sigma}}$ into the column of $\boldsymbol{E}$ corresponding to $v$. 
\item For the rows of type V, we replace the components $\psi(\sigma) - \chi(\sigma)$ by the elements $\boldsymbol{x}_\sigma + \boldsymbol{\nu_{\sigma}}$.
\end{itemize}
To motivate these replacements, note that
\[ \psi(\sigma) - \chi(\sigma) \equiv -\tr(\rho_{\sigma}) = \nu_{\sigma}. \]
It follows that $\pi(\boldsymbol{E}) - E$ has coefficients in $\tilde{I}$.  To see this for the rows of type IV, we use the assumption that $(\xi_v)|_{\cI_v} \equiv \chi|_{\cI_v} \pmod{\tilde{I}}$. For the rows of type V, we use the assumption $\xi_v \equiv \psi|_{\cG_v} \pmod{\tilde{I}}$.
Since $E$ has coefficients in $\T$, the fact that  $\pi(\boldsymbol{E}) - E$ has coefficients in $\tilde{I}$ implies that 
\begin{equation} \label{e:econg}
 \pi(\det(\boldsymbol{E})) \equiv \det(E) \pmod{\tilde{I}}. \end{equation}

We next define a matrix $\boldsymbol{E'}$ from $\boldsymbol{E}$ by enacting the same alterations used to define $E'$ from $E$. More precisely:
\begin{itemize} \item For the first $t$ rows of $\boldsymbol{E}$, the rows of type IV, and the rows of type V, replace
 \[\boldsymbol{x}_\sigma + \boldsymbol{\nu_{\sigma}}  \quad \text{ by } \quad 
 \boldsymbol{x}_\sigma - \boldsymbol{a}_{\sigma}. \]
  \item For rows of type I or III, no change.
 \item For rows of type II,  replace $\boldsymbol{\delta}_{ijj}$ by $\boldsymbol{\delta}_{ijj} - \boldsymbol{a}_i - \boldsymbol{\nu}_i$, replace $\boldsymbol{\delta}_{iji}$ by $\boldsymbol{\delta}_{iji} - \boldsymbol{d}_j$, and leave the other $\boldsymbol{\delta}_{ijk}$ unchanged.
\end{itemize}
We then have $\pi(\boldsymbol{E}') = E'$.  Combining with (\ref{e:econg}), in order to prove our desired result $\det(E') - \det(E) \in \tilde{I}$, it  suffices to prove that $\pi(\det(\boldsymbol{E}') - \det(\boldsymbol{E})) \in \tilde{I}$.

\subsection{Relation Ideal}

We now define the polynomial relations that allow us to reduce $\det(\boldsymbol{E}') - \det(\boldsymbol{E})$ to an expression involving traces and determinants, as in the example of \S\ref{s:example}.  
Define \[ \boldsymbol{\rho_i} = \mat{\boldsymbol{a}_i}{\boldsymbol{b}_i}{\boldsymbol{c}_i}{\boldsymbol{d}_i} \in M_2(R). \]
We let $J \subset R$ be the ideal generated by the following:
\begin{itemize}
\item The 4 coefficients of 
\begin{equation} 
\sum_{i=1}^r \boldsymbol{\epsilon}_i \boldsymbol{\rho}_i \label{e:b1}
\end{equation}
 for each row of type I in $D$.
\item The 4 coefficients of 
\begin{equation} \label{e:b2}
(\boldsymbol{\rho}_i + \boldsymbol{\nu}_i) \boldsymbol{\rho}_j - \sum_{k = 1}^r \boldsymbol{\delta}_{ijk} \boldsymbol{\rho}_k
\end{equation} for each row of type II in $D$.
\

\item For each $v \in \cP$, let $\cB_v$ denote the union of the singleton $\{\sigma_v \}$ with the set of $\sigma$ indexing the rows of type IV in $D$.  For each $v \in \Sigma \setminus \{v_0\}$, let $\cB_v$ denote the set of $\sigma$ indexing the rows of type V in $D$.

For $v \in S\setminus\{v_0\}$ and each $\sigma \in \mathcal{B}_v$, recall the variables  $\boldsymbol{a}_\sigma, \boldsymbol{b}_\sigma, \boldsymbol{c}_\sigma, \boldsymbol{d}_\sigma$ corresponding to the entries of $\boldsymbol{\rho}_{\sigma} (=\boldsymbol{\rho}_i$ for some $i$).
For each pair of distinct $\sigma, \tau \in \mathcal{B}_v$ we include as generators for $J$ the 4 coefficients of
\begin{equation} \label{e:b4}
 \mat{A(\sigma, \tau)}{ B(\sigma, \tau) }{C(\sigma,\tau)}{D(\sigma, \tau)}
\end{equation}
where 
\begin{equation} \label{e:ABCD}
\begin{split}
A(\sigma, \tau) & = \boldsymbol{b}_\sigma\boldsymbol{c}_\tau - (\boldsymbol{x}_{\tau} - \boldsymbol{d}_\tau)(\boldsymbol{x}_{\sigma} - \boldsymbol{a}_\sigma),
 \\
B(\sigma, \tau) & = \boldsymbol{b}_\sigma(\boldsymbol{x}_{\tau} - \boldsymbol{a}_\tau) - 
\boldsymbol{b}_\tau(\boldsymbol{x}_{\sigma} - \boldsymbol{a}_\sigma), \\
C(\sigma, \tau) & = \boldsymbol{c}_\sigma(\boldsymbol{x}_{\tau} - \boldsymbol{d}_\tau) - 
\boldsymbol{c}_\tau(\boldsymbol{x}_{\sigma} - \boldsymbol{d}_\sigma),  \\
D(\sigma, \tau) & = \boldsymbol{b}_\tau\boldsymbol{c}_\sigma - (\boldsymbol{x}_{\sigma} - \boldsymbol{d}_\sigma)(\boldsymbol{x}_{\tau} - \boldsymbol{a}_\tau).
\end{split}
\end{equation}
Of course, the elements $D(\sigma, \tau)$ are redundant since $D(\sigma, \tau) = A(\tau, \sigma)$, but our reasons for describing these generators in terms of the matrix (\ref{e:b4}) will become apparent when we prove Lemma~\ref{l:stable} below.
\end{itemize}

\begin{lemma} We have $\pi(J) =0$. \end{lemma}

\begin{proof}
The fact that the elements in (\ref{e:b1}) and (\ref{e:b2}) map to 0 under $\pi$ is clear since these are precisely the relations defining $\epsilon_i, \delta_{ijk}$ in $\T$.  For the generators in (\ref{e:b4}) we compute as follows.  
\begin{equation} \label{e:pibst}
 \pi(B(\sigma, \tau)) = b(\sigma)(\xi_v(\tau) - \psi(\tau) - a(\tau)) - b(\tau)(\xi_v(\sigma) - \psi(\sigma) - a(\sigma)).
\end{equation}
Using (\ref{e:db}), it is easy to see that the expression (\ref{e:pibst}) vanishes if it is multiplied by $B_v$ or by $D_v$.
Since $B_v$ and $D_v$ generate $K$ as a $K$-module, it follows that (\ref{e:pibst}) vanishes.  The proof that $A(\sigma, \tau), C(\sigma, \tau),$ and $D(\sigma, \tau)$ also vanish under $\pi$ is similar and uses the equation
\[
 B_v c(\sigma) = D_v(\xi_v(\sigma) - \psi(\sigma) - d(\sigma))
\]
obtained from equating the lower right entries of (\ref{e:changebasis}).
\end{proof}

\subsection{Subring of traces and determinants} \label{s:subring}

In order to apply the results on traces and determinants from \S\ref{s:ti}, we define a certain subring $A \subset R$ as follows.
Let $A_1 \subset R_1$ denote the sub $R_0$-algebra generated by the traces and determinants of all matrices in the noncommutative $\Z$-algebra generated by the matrices $\boldsymbol{\rho}_i$:
\[ A_1 = R_0[\boldsymbol{a}_i + \boldsymbol{d}_i, \boldsymbol{a}_i\boldsymbol{d}_i - \boldsymbol{b}_i\boldsymbol{c}_i, 
\boldsymbol{a}_i\boldsymbol{a}_j + \boldsymbol{b}_i\boldsymbol{c}_j + \boldsymbol{c}_i\boldsymbol{b}_j + \boldsymbol{d}_i\boldsymbol{d}_j,  \dotsc ]. \]
Let $A_0$ denote the image of $A_1$ in $R = R_1/(\boldsymbol{b}_{\tau} \colon \tau \in \cB_{v_0})$. 
Let $A \subset R$ denote the ring generated over $A_0$ by the $\boldsymbol{d}_\tau$ for $\tau \in \cB_{v_0}$.
Denote by $\overline{A}$ the image of $A$ in $R/J$.  

The goal of the remainder of this subsection is to show that in order to deduce our desired result  $\pi(\det(\boldsymbol{E}') - \det(\boldsymbol{E})) \in \tilde{I}$, it suffices to prove that the image of $\det(\boldsymbol{E}') - \det(\boldsymbol{E})$ in $R/J$ lies in $\overline{A}$.  For this, it is important that $\det(\boldsymbol{E}') - \det(\boldsymbol{E})$ lies in the following ideal of $R$:
\[ I_R = \langle \boldsymbol{a}_i + \boldsymbol{\nu}_i,   \boldsymbol{b}_i,  \boldsymbol{c}_i,  \boldsymbol{d}_i   \rangle. \]
The fact that  $\det(\boldsymbol{E}') - \det(\boldsymbol{E})$ lies in $I_R$ follows from multilinearity of the determinant since every entry of 
$\boldsymbol{E}' - \boldsymbol{E}$ lies in $I_R$.

\begin{lemma} We have $\pi(A \cap (I_R, J)) \subset \tilde{I}$. \label{l:pia}
\end{lemma}

\begin{proof}
Let  $\cB_{v_0} = \{\tau_1, \dots, \tau_h\}$. 
Define a $R_0$-algebra homomorphism
\[ \boldsymbol{V} \colon R_0 \langle X_1, \dots, X_r \rangle \longrightarrow R_0, \qquad  X_i \mapsto -\boldsymbol{\nu_i}. \]
Any element of $A$ can be written as a polynomial in the expressions
\begin{align}
&  \tr(f(\boldsymbol{\rho}_1, \dots, \boldsymbol{\rho}_r)) - \boldsymbol{V}(f), \label{e:trv} \\
&  \det(f(\boldsymbol{\rho}_1, \dots, \boldsymbol{\rho}_r)), \label{e:detf} \\
& \text{ and } \boldsymbol{d}_{\tau_1}, \dots, \boldsymbol{d}_{\tau_h} \label{e:ds} \end{align}
  with coefficients in $R_0$, as $f$ ranges over polynomials with no constant term.  As we now explain, the elements in (\ref{e:trv})--(\ref{e:ds}) have images under $\pi$ lying in $\tilde{I}$.  This holds by  Lemma \ref{l:tr-char} for (\ref{e:trv}) and 
  by Lemma~\ref{l:dets} for (\ref{e:detf}).  For (\ref{e:ds})
 we note that our choice of basis ensures that
  \[ \pi(\boldsymbol{d}_{\tau_i}) = \xi_{v_0}(\tau_i) - \psi(\tau_i) \]
(see \ref{e:rholocal}), and by assumption (\ref{e:Sigmaassumption}) this element lies in $\tilde{I}$.
  It is also easy to check that each  expression (\ref{e:trv})--(\ref{e:ds}) lies in $I_R$.  For example, we note that modulo $I_R$ we have $\boldsymbol{\rho}_i \equiv \mat{-\boldsymbol{\nu}_i}{0}{0}{0}$,
  from which it follows that  (\ref{e:trv}) lies in $I_R$. The computation that (\ref{e:detf})  lies in $I_R$ is similar.  The elements in (\ref{e:ds}) lie in $I_R$ by definition.
  
  To conclude, it therefore suffices to prove that $\pi(R_0 \cap (I_R, J)) \subset \tilde{I}$.  To calculate the intersection $R_0 \cap (I_R, J)$, we note that
  \[ R / I_R \cong R_0, \] 
  with $\boldsymbol{a}_i \mapsto - \boldsymbol{\nu_i}$ and $\boldsymbol{b}_i, \boldsymbol{c}_i, \boldsymbol{d}_i \mapsto
  0$.  
  To compute the quotient by $J$ we then mod out by the elements of $J$ with these substitutions made, and obtain
  \[ R/(I_R, J) \cong R_0/I_{R_0}
  \]
  where $I_{R_0} \subset R_0$ is the ideal generated by
  \begin{itemize}
 \item  $\sum_{i=1}^r \boldsymbol{\epsilon}_i \boldsymbol{\nu}_i$ for each row of type I,
 \item $\sum_{k=1}^r \boldsymbol{\delta}_{ijk} \boldsymbol{\nu}_k$ for each row of type II,
\item $(\boldsymbol{x}_\sigma + \boldsymbol{\nu}_{\sigma}) \boldsymbol{x}_\tau$ for each distinct $\sigma, \tau \in \cB_v$, for each $v \in S \setminus \{v_0\}.$
\end{itemize}
Thus $R_0 \cap (I_R, J) = I_{R_0}$ and it remains to check that each of the bulleted expressions maps to $\tilde{I}$ under $\pi$.  For instance, for a row of type I we have
\[ 0 = \tr\Big(\sum_{i=1}^r  {\epsilon}_i {\rho}_i \Big)
\equiv - \sum_{i=1}^r {\epsilon}_i {\nu}_i \pmod{I} \]
by Lemma~\ref{l:tr-char}.  The proof for the second bulleted item is similar.  For the third, we note that 
\begin{equation} \label{e:pisc}
 \pi((\boldsymbol{x}_\sigma + \boldsymbol{\nu}_{\sigma}) \boldsymbol{x}_\tau) = (\xi_v(\sigma) - \chi(\sigma))(\xi_v(\psi) - \psi(\sigma)). 
 \end{equation}
All terms in (\ref{e:pisc}) lie in $\tilde{\T}$.  Furthermore $\xi_v(\sigma) - \chi(\sigma) \in \tilde{I}$ for $v \in \cP$ and $\xi_v(\sigma) - \psi(\sigma) \in \tilde{I}$ for $v \in \Sigma$ by assumptions (\ref{e:Passumption}) and (\ref{e:Sigmaassumption}).  This concludes the proof.
\end{proof}

We summarize the result of this section.

\begin{prop}
 In order to prove the desired result 
 \[ \det(E') - \det(E) \in \tilde{I}, \]
 it suffices to prove that the image of 
$e = \det(\boldsymbol{E}') - \det(\boldsymbol{E})$ in $R/J$ lies in the subring $\overline{A}$, the image of $A$ in $R/J$.
\end{prop}
\begin{proof} If the image of $e$ in $R/J$ lies in $\overline{A}$, then we may write $e = a+ j$ with $a \in A$ and $j \in J$.  Then $a = e - j \in A \cap (I_R, J)$ since $e \in I_R$, so Lemma~\ref{l:pia} implies that $\pi(a) \in \tilde{I}$.  Since $\pi(j) = 0$, we obtain $\pi(e) \in \tilde{I}$ as well.  We showed in \S\ref{s:trr} that $\pi(e) \equiv \det(E') - \det(E) \pmod{\tilde{I}}$.  The result follows.
\end{proof}

\section{Matrix Invariant Theory and Rational Cohomology}

As mentioned above, we do not know how to use direct computation to show that the image of $\det(\boldsymbol{E}') - \det(\boldsymbol{E})$ in $R/J$ lies in $\overline{A}$.  Instead, we identify the subring $A \subset R$ as the ring of invariants of a Borel subgroup of the algebraic group 
$\GL_2\!/\Z$ acting by simultaneous conjugation on the matrices $\boldsymbol{\rho}_i$.  The cohomology theory of algebraic group actions goes by the name {\em rational cohomology}, or alternately {\em Hochschild cohomology}, and for the convenience of the reader we review some basic definitions.  As mentioned earlier, the word ``rational'' is used in the sense of rational maps, not the rational numbers; we work integrally over $\Z$ throughout.

\subsection{Rational cohomology}
Let $\G$ denote the algebraic group $\GL_2\!/\Z$, i.e.
\[ \G = \Spec \mathfrak{G}, \qquad \mathfrak{G} = \Z[a, b, c, d, (ad-bc)^{-1}] \]
endowed with its usual structure as a group scheme given by a comultiplication
\[ c\colon \mathfrak{G} \longrightarrow \mathfrak{G} \otimes_{\Z} \mathfrak{G} \]
and a counit $e \colon \mathfrak{G} \lra \Z$.
 The lower triangular Borel $\B \subset \G$ will play an important role in this study:
\[ \B = \Spec \mathfrak{B}, \qquad \mathfrak{B} = \Z[x, y, z, (xz)^{-1}], \]
with $\mathfrak{G} \rightarrow \mathfrak{B}$ given by $\mat{a}{b}{c}{d} \mapsto \mat{x}{0}{y}{z}$.

\begin{definition}[\cite{jantzen}*{I.2.7-I.2.8}]  A {\em $\G$-module} is a $\Z$-module $V$ endowed with a \emph{comodule map} $w \colon V \longrightarrow \mathfrak{G} \otimes_{\Z} V$, a map of $\Z$-modules satisfying:
\begin{itemize}
\item $(1 \otimes w) \circ w = (c \otimes 1) \circ w$,
\item $(e \otimes 1) \circ w = id$.
\end{itemize}

Equivalently, a $\G$-module is a $\Z$-module $V$ endowed with a functorial action, for each commutative ring $A$, of $\G(A) = \GL_2(A)$ on $A \otimes_{\Z} V$. Consider the universal element \[ g =\mat{a}{b}{c}{d} \in \G(\mathfrak{G}),\] corresponding to $id$ under the isomorphism $\G(\mathfrak{G}) \cong \Hom(\mathfrak{G}, \mathfrak{G})$. The action of $g$ on $\mathfrak{G} \otimes_{\Z} V$ defines a map $\mathfrak{G} \otimes_{\Z} V  \lra \mathfrak{G} \otimes_{\Z} V$ whose restriction to $1 \otimes V$ is $w$. 
\

  We say that $V$ is a $\G$-module over a commutative ring $R$ if $V$ has a structure of $R$-module that commutes with the $\G$-action.
\end{definition}

The center $\G_m \subset \G$ decomposes any $\G$-module $V$ into a sum of \emph{homogeneous} $\G$-modules: $V = \oplus_{n \in \Z} V^{(n)}$, where $V^{(n)} = \{v \in V\colon \mat{x}{0}{0}{x} v = x^n \otimes v\}$.

\bigskip

\begin{example} \label{e:adjoint}
 We denote by $\cA = \Z A \oplus \Z B \oplus \Z C \oplus \Z D$ the $\G$-module given by the adjoint representation, i.e. for $g = \mat{a}{b}{c}{d}$ the universal element,
\[ \mat{g \cdot A}{g \cdot B}{g \cdot C}{g \cdot D} = \mat{a}{b}{c}{d}^{-1} \mat{A}{B}{C}{D}  \mat{a}{b}{c}{d}. \]

Note that for $g = \mat{x}{0}{y}{z} \in \B$, we have \[ g \cdot A = A + \frac{y}{x} B, \quad g \cdot B = \frac{z}{x} B. \] In particular, \[ \mathcal{B} = \Z B \quad \text{ and }  \quad \mathcal{V} =  \Z A \oplus \Z B \] are $\B$-submodules of $\cA$.
For future reference, we also record \[ g \cdot C = -\frac{y}{z}A - \frac{y^2}{xz}B + \frac{x}{z}C + \frac{y}{z} D, \qquad g \cdot D = -\frac{y}{x} B + D. \]

For a $\B$-module $M$, let 
\begin{equation} \label{e:twistdef}
M(n) := M \otimes \mathcal{B}^{\otimes n}. \end{equation} Note that $\bigwedge^2 \mathcal{V} \cong \Z(1)$ and $\mathcal{V}^* \cong \mathcal{V}(-1)$. 

\end{example}

\begin{definition}[\cite{jantzen}*{I. 4}]  If $V$ is a $\G$-module, define the invariants
\begin{align*} H^0(\G, V) = V^{\G} &= \{v \in V \colon v \mapsto 1 \otimes v \text{ under the coaction map } V \rightarrow \mathfrak{G} \otimes V \} \\
&= \{ v \in V \colon \text{for all commutative rings $A$ and $g \in \G(A)$, } g \cdot v = v\}.  \end{align*}
The rational cohomology groups $H^i(\G, -)$ are the right derived functors of $H^0(\G, -)$.  \end{definition}

All of the definitions above carry through with $\G$ replaced by $\B$.
The following is proven in \cite{cpsv}*{Theorem 2.1} for semisimple groups over algebraically closed fields in characteristic $p$. However, these hypotheses can be removed using the techniques of \cite{jantzen}.
\begin{theorem} \label{t:bg} Let $V$ be a $\G$-module.  The restriction map \[ H^i(\G, V) \lra H^i(\B, V) \] is an isomorphism for all $i \ge 0$. 
\end{theorem}

\begin{proof} 
For a field $k$, the isomorphism $H^i(\G_k, V_k) \cong H^i(\B_k, V_k)$ is \cite{jantzen}*{II.4.7}. To deduce the result over $\Z$, apply the universal coefficient theorem \cite{jantzen}*{I.4.18} (see also the proof of \cite{jantzen}*{II.4.5}).
\end{proof}

\begin{remark} Theorem~\ref{t:bg} is our motivation for using rational cohomology.  The analogous statement in ordinary group cohomology for the groups $\G(\Z) = \GL_2(\Z)$ and $\B(\Z)$ is in general false.
\end{remark}

\subsection{Acyclic $\B$-modules}

Given a $\B$-module $W$, the algebraic induction $\Ind_{\B}^{\G}(W)$ is defined to be the $\G$-module $(\mathfrak{G} \otimes_{\Z} W)^{\B}$, where the action of $\B$ on $\mathfrak{G}$ comes from right-multiplication on $\G$ (\cite{jantzen} I.3.3).

\begin{definition} Let $V$ be a homogeneous $\G$-module.  We say that $V$ has a \emph{good filtration} if there is a finite filtration $0 = V_0 \subset \cdots \subset V_n = V$ such each $V_{i+1}/V_i$ is isomorphic to the algebraic induction $\Ind_{\B}^{\G}(- \lambda)$ for a dominant weight $\lambda$ (considered as a rank one $\B$-module). More generally, a $\G$-module $V = \oplus_{n \in \Z} V^{(n)}$ is said to have a good filtration if each homogeneous component $V^{(n)}$ has a good filtration.
 \end{definition}
 
 The following result is proven in \cite{cpsv} in characteristic $p$. The proof easily generalizes to $\G = \GL_2\!/\Z$ using results of \cite{jantzen}.
\begin{theorem}\label{t:good-acyclic}
$\G$-modules $V$ with a good filtration are acyclic: $H^i(\G, V) = 0$ for $i > 0$.
\end{theorem} 
\begin{proof}
It suffices to prove for $V =\Ind_{\B}^{\G}(- \lambda)$. The right derived functors of induction, $R^*\Ind_{\B}^{\G}(-\lambda)$, are isomorphic to sheaf cohomology $H^*(\G/\B, \mathcal{L}(-\lambda))$ of a certain coherent sheaf $\mathcal{L}(-\lambda)$ (\cite{jantzen}*{I.5.13}) . In fact, $\mathcal{L}(-\lambda) \cong \cO(n)$ for some $n > 0$, hence the higher cohomology vanishes: $R^i\Ind_{\B}^{\G}(-\lambda) \cong H^i(\mathbf{P}^1, \cO(n)) = 0$ for $i > 0$. Therefore one can apply \cite{jantzen}*{I.4.6} to obtain $H^*(\G,\Ind_{\B}^{\G}(- \lambda)) = H^*(\B, -\lambda)$. This equals 0, since the weight $-\lambda$ is not a positive root (see \cite{cpsv}*{2.2}).
\end{proof}

The following result is well-known (see \cite{cp}*{Theorem 11.1}).
\begin{theorem}\label{t:tensor-good}
If $V$ and $W$ are $\G$-modules with good filtrations, so is $V \otimes_{\Z} W$.
\end{theorem} 
 
\begin{corollary}\label{c:adj-good}
The adjoint representation $\cA$ is a good $\G$-module, as are its tensor powers $\cA^{\otimes k}$.
\end{corollary}
\begin{proof}
We confirm that $\cA$ has a good filtration. The standard representation $\mathrm{std}$ is induced from the dominant weight $\lambda = (0, -1)$, while the dual of the standard $\mathrm{std}^{\vee}$ is induced from the dominant weight $\lambda =(1,0)$. Thus $\cA \cong \mathrm{std} \otimes \mathrm{std}^{\vee}$ has a good filtration.
\end{proof}

We also need the following result (\cite{cp} Corollary 11.7):
\begin{theorem}\label{t:r-good}
The module $\mathbf{Z}[\boldsymbol{a}_i,  \boldsymbol{b}_i, \boldsymbol{c}_i, \boldsymbol{d}_i]_{i=1,\ldots,n}$, where $\G$ acts via simultaneous conjugation, has a good filtration.
\end{theorem}

\begin{lemma} \label{l:vanishing}
For any good $\G$-module $V$ and $j \geq 0$, $H^i(\B, V(j)) = 0$ for $i > j$.
\end{lemma}
\begin{proof}
If $j = 0$, this follows from Theorem \ref{t:good-acyclic} and Theorem \ref{t:bg}. If $j \geq 1$, we have \[ H^i(\B, V(j)) \cong H^{i-1}(\G, V\otimes H^1(\mathbf{P}^1, \cO(-2j))) \cong H^{i-1}(\G, V\otimes \Sym^{2j-2}(\mathrm{std})^{\vee})\] (see \cite{jantzen}*{I.4.5 and I.4.8}). We thus need to show that $H^{i-1}(\G, V \otimes \Sym^{2j-2}(\mathrm{std})^{\vee}) = 0$ when $i > j$. By \cite{jantzen}*{I.4.18}, it suffices to prove that for all primes $p > 0$, we have \[ H^{i-1}(\GL_{2}\!/{\F_p}, V \otimes \Sym^{2j-2}(\mathrm{std})^{\vee} \otimes_{\Z} \F_p) = 0\]  when $i > j$. 
By \cite{fp}*{\S3}, this is equivalent to proving that the \emph{good filtration dimension} of \[ \Sym^{2j-2}(\mathrm{std})^{\vee} \otimes \F_p \] is $\leq j-1$ for all primes $p$. By \cite{p}*{Theorem 4.2 and Lemma 3.5}, $\Sym^{2j-2}(\mathrm{std})^{\vee} \otimes \F_p$ has good filtration dimension equal to $\lfloor \frac{2j-2}{p} \rfloor$.
\end{proof}

\begin{lemma}\label{l:higher-cohom}
For any good $\G$-module $V$,
\begin{enumerate}
\item $H^1(\B, V(1)) \cong H^0(\G,  V)$,
\item for $i > 1$, $H^i(\B,  V(i)) \otimes \F_p \cong 0$ for primes $p > 2$,
\item for $i > 1$,  $H^i(\B,  V(i)) \otimes \F_2 \cong H^0(\G, V) \otimes \F_2$.
\end{enumerate} 
\end{lemma}
\begin{proof}
As in Lemma \ref{l:vanishing}, we have an isomorphism \[ H^i(\B, V(i))  \cong H^{i-1}(\G, V\otimes \Sym^{2i-2}(\mathrm{std})^{\vee}). \]  Taking $i = 1$ proves (1). 

For a prime $p > 2$,  $H^{i-1}(\G, V\otimes \Sym^{2i-2}(\mathrm{std})^{\vee}) \otimes \F_p = 0$, as $i-1$ is greater than the good filtration dimension of $\Sym^{2ji-2}(\mathrm{std})^{\vee}$.  This proves (2).

For $p = 2$,
\[ H^{i-1}(\G, V\otimes \Sym^{2i-2}(\mathrm{std})^{\vee}) \otimes \F_2 \cong \Ext^{i-1}_{\G}(\Sym^{2i-2}(\mathrm{std}) \otimes \F_2, V \otimes \F_2). \] By repeatedly applying the final isomorphism in the proof of \cite{p}*{Theorem 4.2},  we find that \[ \Ext^{i-1}(\Sym^{2i-2}(\mathrm{std}) \otimes \F_2, \Ind_\B^{\G}(-\nu) \otimes \F_2) \cong \Ext^0(\F_2, \Ind_\B^{\G}(-\nu) \otimes \F_2), \] where $\nu$ is a dominant weight. As $V$ admits a filtration whose gradeds are such $\Ind_\B^{\G}(-\nu)$,  this implies that
\[  \Ext^{i-1}(\Sym^{2i-2}(\mathrm{std}) \otimes \F_2,V \otimes \F_2) \cong  \Ext^{0}(\F_2,V \otimes \F_2) = H^0(\G,  V \otimes \F_2). \]

Using the long exact sequence associated to $0 \to V \xrightarrow{\cdot 2} V \to V \otimes \F_2 \to 0$, combined with $H^1(\G,V) = 0$ (Lemma \ref{l:vanishing}), we find that $H^0(\G, V \otimes \F_2) \cong H^0(\G, V) \otimes \F_2$. This concludes the proof of (3).
\end{proof}

\begin{corollary}\label{c:base-change}
Given an inclusion $S_1 \subset S_2$ of commutative rings which are good $\G$-modules, the natural map \[ H^i(\B, S_1(i)) \otimes_{H^0(\B,S_1)} H^0(\B, S_2) \to H^i(\B, S_2(i)) \] is surjective for all $i \geq 0$.
\end{corollary}
\begin{proof}
The case $i = 0$ is automatic, and the case $i = 1$ follows from Lemma \ref{l:higher-cohom} (1). 

For $i > 2$,  we must show that this map is surjective after tensoring with $\F_p$ for any prime $p$. By Lemma \ref{l:higher-cohom} (2),  the domain and codomain both vanish unless $p = 2$.

The map we are considering is defined via the product structure in group cohomology. As all of the degree-reducing isomorphisms in Lemma \ref{l:higher-cohom} came from (the inverse of) boundary maps associated to certain long exact sequences, they are compatible with this product structure. Thus it suffices to check that the natural map
\[ (H^0(\G, S_1) \otimes \F_2) \otimes_{H^0(\G,S_1)} H^0(\G, S_2) \to (H^0(\G, S_2) \otimes \F_2) \] is surjective.  This is clear.
\end{proof}

\begin{lemma}\label{l:nice-res}
If a $\B$-module $V$ admits a finite resolution \[ \begin{tikzcd}
 V_n(n) \ar[r,"f_n"] & V_{n-1}(n-1)  \ar[r,"f_{n-1}"] & \cdots \ar[r, "f_2"] & V_1 \ar[r,"f_1"] & V_0 \ar[r,"f_0"] & V
 \end{tikzcd} \] where each $V_j$ is a good $\G$-module, then $V$ is an acyclic $\B$-module.
\end{lemma}
\begin{proof}
For $j=0, \dots, n+1$, write $W_j = \ker(f_{j-1}) = \im(f_{j})$, so in particular $W_{0} = V$ and $W_{n+1} =0$.  For $j=0, \dots, n$, we have a short exact sequence
\[ 0 \lra W_{j+1} \lra V_{j}(j) \lra W_{j} \lra 0. \]
Taking the associated long exact sequences in $\B$-cohomology, a downward induction using Lemma~\ref{l:vanishing} shows that $H^i(\B, W_j)$ vanishes for $i > j$. 
For $j=0$, this is the desired result.
\end{proof}

\begin{theorem}\label{t:r-acyclic}
Consider the $\B$-module \[ S = \Z[\boldsymbol{a}_i, \boldsymbol{b}_i, \boldsymbol{c}_i, \boldsymbol{d}_i]_{i=1,\ldots,n}/(\boldsymbol{b}_i \colon 1 \leq i \leq k), \] where $\B$ acts by simultaneous conjugation. Let $V$ be a good $\G$-module such that $V$ is $\Z$-flat.  Then $V \otimes_{\Z} S$ is an acyclic $\B$-module.
\end{theorem}
\begin{proof}
Let $S_0 = \Z[\boldsymbol{a}_i, \boldsymbol{b}_i, \boldsymbol{c}_i, \boldsymbol{d}_i]_{i=1,\ldots,n}$. The elements $\boldsymbol{b}_1,\ldots, \boldsymbol{b}_k$ form a regular sequence in $S_0$, hence the Koszul complex for this sequence yields an exact sequence of $\B$-modules:
\[ \left(\bigwedge^k_{S_0} \bigoplus_{i=1}^k S_0\right)(k) \lra \cdots \lra \left(\bigoplus_{i=1}^k S_0\right)(1) \lra S_0. \]
This gives a resolution of $S = S_0/(\boldsymbol{b}_1, \dotsc, \boldsymbol{b}_k)$. If we tensor this complex with $V$, we obtain an exact complex of $\B$-modules resolving $V \otimes_\Z S$. The $k$th term of this complex is of the form $\left(\bigoplus V \otimes S_0 \right)(k)$. As $V \otimes S_0$ is a good $\G$-module by Theorem \ref{t:r-good} and Theorem \ref{t:tensor-good}, we may apply Lemma \ref{l:nice-res} to conclude that $V \otimes S$ is an acyclic $\B$-module.
\end{proof}

In particular, this theorem may be applied to the ring $R$ defined in \S\ref{s:formal}, with the action of $\B$ defined in \S\ref{s:invt-theory}.

\subsection{Matrix Invariant Theory}\label{s:invt-theory}

The following important classical result is the subject of the beautiful book \cite{cp} by de Concini and Procesi.  

\begin{theorem}\label{t:invariants}Let $R_0$ be a $\Z$-flat commutative ring and let \[ R = R_0[a_i, b_i, c_i, d_i]_{i=1}^{r}. \] Endow $R$ with a $\G$-action defined by simultaneous conjugation on the matrices $\rho_i = \mat{a_i}{b_i}{c_i}{d_i}$. The ring $H^0(\G, R) \subset R$ is generated as an algebra over $R_0$ by the traces and determinants of all matrices in the algebra generated by the $\rho_i$.
\end{theorem}

Theorem~1.10 of {\em loc.\ cit.}\ is the statement above with $R_0 = \Z$ and invariants taken for the group $\G(\Z) = \GL_2(\Z)$.  Since $\G$-invariance is stronger than $\G(\Z)$-invariance, the statement above follows for $R_0 = \Z$.  Extension to arbitrary $\Z$-flat algebras $R_0$ on which $\G$ acts trivially is then immediate.

We now apply Theorem~\ref{t:invariants}  to our setting.
Recall the ring $R$ defined in (\ref{e:r0def})--(\ref{e:rdef}) above and the subring $A \subset R$ defined at the start of \S\ref{s:subring}.
Let $\B$ act on $R_1 = R_0[\boldsymbol{a}_i, \boldsymbol{b}_i, \boldsymbol{c}_i, \boldsymbol{d}_i]_{i=1}^r $ by simultaneous conjugation on the matrices $\boldsymbol{\rho}_i$. In particular, $\B$ acts trivially on $R_0.$ This action preserves the ideal $(\boldsymbol{b}_{\tau} \colon \tau \in \mathcal{B}_{v_0})$ (see Example \ref{e:adjoint}), and hence descends to an action on $R = R_1/(\boldsymbol{b}_{\tau} \colon \tau \in \mathcal{B}_{v_0}).$

\begin{corollary}\label{c:a} 
We have $H^0(\B, R) = A$.
\end{corollary}

\begin{proof}
It is convenient to write $R$  in the form
\begin{equation} \label{e:taunotation}
 R = R_0[\boldsymbol{a}_i, \boldsymbol{b}_i, \boldsymbol{c}_i, \boldsymbol{d}_i]_{i = 1}^s[\boldsymbol{a_{\tau}},\boldsymbol{b_{\tau}}, \boldsymbol{c_{\tau}}, \boldsymbol{d_{\tau}}]/(\boldsymbol{b}_{\tau}) \end{equation}
where we exclude the variables corresponding to the entries of $\boldsymbol{\rho}_{\tau}$ with $\tau \in \mathcal{B}_{v_0}$ from the first list of variables, and $\tau$ ranges over $\mathcal{B}_{v_0}$ in the second set of variables.
We first handle the case $s=0$, i.e.\ 
\[ R = R_0[\boldsymbol{a_{\tau}},\boldsymbol{b_{\tau}}, \boldsymbol{c_{\tau}}, \boldsymbol{d_{\tau}}]/(\boldsymbol{b}_{\tau}). \]
We claim that in this case $H^0(\B, R)$ is the image of $R_0[\boldsymbol{a_{\tau}},\boldsymbol{d_{\tau}}]$ in $R$.  Let $\D \subset \B$ denote the torus of diagonal matrices.  There is a $\D$-module isomorphism
\[ R \cong R_0[\boldsymbol{a_{\tau}}, \boldsymbol{c_{\tau}}, \boldsymbol{d_{\tau}}], \qquad  \boldsymbol{b_{\tau}} \mapsto 0. \]
Since $\D$ acts trivially on $\boldsymbol{a_{\tau}}$ and $\boldsymbol{d_{\tau}}$ but scales $\boldsymbol{c_{\tau}}$ non-trivially, it is clear that 
$H^0(\D, R)$ is the image of  $R_0[\boldsymbol{a_{\tau}},\boldsymbol{d_{\tau}}]$.  Since the elements $\boldsymbol{a_{\tau}}$ and $\boldsymbol{d_{\tau}}$ are $\B$-invariant in $R$, the claim follows.

Now we return to the case of general $s$.
Let \[ S = R_0[\boldsymbol{a}_i, \boldsymbol{b}_i, \boldsymbol{c}_i, \boldsymbol{d}_i],  \text{ so that } R_1 = S[\boldsymbol{a_{\tau}},\boldsymbol{b_{\tau}}, \boldsymbol{c_{\tau}}, \boldsymbol{d_{\tau}}]. \] Let $\fb = (\boldsymbol{b}_{\tau} \colon \tau \in \mathcal{B}_{v_0})$.  By Theorem \ref{t:invariants},  the image of $H^0(\G, R_1) = H^0(\B,R_1)$ in $H^0(\B, R_1/\fb)$ is the ring $A_0$ defined in \S\ref{s:subring}.  

Let $R'_1 =R_0[\boldsymbol{a_{\tau}},\boldsymbol{b_{\tau}}, \boldsymbol{c_{\tau}}, \boldsymbol{d_{\tau}}] \subset R_1$, and let $\fb' = R'_1 \cap \fb. $ By the claim above for $s=0$, we have \[ H^0(\B, R'_1/\fb') = R_0[\boldsymbol{a_{\tau}},\boldsymbol{d_{\tau}}]. \]
 The ring $A$ is by definition generated over $A_0$ by the $\boldsymbol{d_{\tau}}$'s and hence also by the larger generating set of $\boldsymbol{a_{\tau}}$'s and $\boldsymbol{d_{\tau}}$'s (note $\boldsymbol{a_{\tau}} + \boldsymbol{d_{\tau}} \in A_0$). Thus, to prove that $H^0(\B, R_1/\fb) = A$,  it will suffice to show that the map \begin{equation}\label{e:base-change} H^0(\B, R'_1/\fb) \otimes_{H^0(\B, R'_1)} H^0(\B, R_1) \lra H^0(\B, R_1/\fb)\end{equation} is surjective.

 Since  $\{\boldsymbol{b}_\tau \}$ is a regular sequence in $R_1$, we may provide a $\B$-module resolution of $R_1/\fb$ via a Koszul complex:
\[ \left(\bigwedge\nolimits^k_{R_1}  \left(\bigoplus_{i=1}^k R_1\right)\right)(k) \lra \cdots \lra  \bigoplus_{i=1}^k R_1(1) \lra R_1 \lra R_1/\fb \lra 0.  \]

Here $k$ denotes the number of $\tau$'s, and the map $\bigoplus_{i=1}^k R_1(1) \lra R_1$ sends each standard basis vector on the left to a $\boldsymbol{b_{\tau}}$.
Let \[ C_j =  \left(\bigwedge\nolimits^j_{R_1}  \left(\bigoplus_{i=1}^j R_1\right)\right)(j) \]
denote the terms of this complex.
 There is an analogous Koszul complex for $R'_1/\fb'$, which we denote $C'_j$. We have $C_{\bullet} \cong C'_{\bullet} \otimes_{R'_1} R_1$.

  There is a 2nd-quadrant spectral sequence \[ E_1^{-i,j} = H^j(\B,  C_i) \implies H^{j-i}(\B, R_1/\fb). \]  Since $R_1$ is a good $\G$-module (Theorem \ref{t:r-good}),  Lemma \ref{l:vanishing} yields $H^j(\B, C_i) = 0$ for $j > i$. This vanishing implies that we have a sequence of maps
\begin{align*}
 \alpha_0 &\colon H^0(\B, C_0) \lra H^0(\B, R_1/\fb),  \\
 \alpha_i &\colon H^i(\B,  C_i) \lra H^0(\B, R_1/\fb)/\alpha_{i-1}(H^{i-1}(\B, C_{i-1}))\text{ for } i = 1,\ldots, k,
 \end{align*}
  such that $\alpha_k$ is surjective.  In other words, $\alpha_i(H^i(\B,  C_i))$ are the associated gradeds of the filtration on $H^0(\B, R_1/\fb)$ induced by the spectral sequence. 
We similarly have that $\alpha_i(H^i(\B,  C'_i))$, $i = 0,\ldots, k$, are the associated gradeds of the filtration on $H^0(\B, R'_1/\fb')$ induced by the corresponding spectral sequence for $C'_{\bullet}$.

To prove the surjectivity of (\ref{e:base-change}),  we will show that the natural map \begin{equation} \label{e:alphamap}
 \alpha_i(H^i(\B,  C'_i)) \otimes_{H^0(\B, R'_1)} H^0(\B, R_1) \lra \alpha_i(H^i(\B,  C_i)) \end{equation} is surjective for all $i = 0, \ldots, k$. Since $C'_i \cong \bigoplus R'_1(i)$ and $C_i \cong \bigoplus R_1(i)$,  and both $R'_1$ and $R_1$ are good $\G$-modules, we apply Corollary \ref{c:base-change} to see that the natural map \[ H^{i}(\B, C'_{i}) \otimes H^0(\B, R_1)  \lra H^{i}(\B, C_{i}) \] is surjective. The surjectivity of (\ref{e:alphamap}) follows.

\end{proof}

\subsection{Roadmap}

We now return to the proof of our main result, Theorem~\ref{t:rll}, picking up the argument of \S\ref{s:formal}.
Recall that we had defined an element
\[ \boldsymbol{e} = \det(\boldsymbol{E}') - \det(\boldsymbol{E}) \] and our goal is to show that its
 its image in $R/J$, denoted $\overline{\boldsymbol{e}}$, lies in  $\overline{A} \subset R/J$. Here $\overline{A}$ is the image of the subring $A \subset R$ in $R/J$.
 Our strategy is as follows.  We will first prove that the ideal $J$ is stable under the action of $\B$.
We will then show that 
\[ \overline{\boldsymbol{e}} \in H^0(\B, R/J). \]
The long exact sequence in rational cohomology associated to 
\[ 0 \longrightarrow J \longrightarrow R \longrightarrow R/J \longrightarrow 0 \]
 yields an exact sequence
\begin{equation} \label{e:les}
 \begin{tikzcd}
H^0(\B, R) = A \ar[r] & H^0(\B, R/J) \ar[r, "c_J"] & H^1(\B, J).  \end{tikzcd}
\end{equation} 
The equality on the left is Corollary~\ref{c:a}.  Let $\beta \in H^1(\B, J)$ denote the image of $\overline{\boldsymbol{e}}$ under the connecting homomorphism $c_J$.
In view of (\ref{e:les}), the desired result $\overline{\boldsymbol{e}} \in \overline{A}$ will follow if we can show that $\beta =0$.

For this, we will define a certain $\B$-stable subideal $J' \subset J$ and show that in fact 
\[ \overline{\boldsymbol{e}}\in H^0(\B, R/J'). \]
(This is a slight abuse of notation, as here $\overline{\boldsymbol{e}}$ denotes the reduction of $\boldsymbol{e}$ modulo $J'$.)

We let $\alpha = c_{J'}(\overline{\boldsymbol{e}}) \in H^1(\B, J')$ defined as above.  Then $\beta = \iota_*(\alpha)$ where $\iota \colon J' \rightarrow J$ is the inclusion and $\iota_*$ is the induced map on rational cohomology.  To conclude, we will prove that the map
\[ \iota_*\colon H^1(\B, J') \longrightarrow H^1(\B, J) \]
vanishes.  Therefore $\beta = 0$, and our result follows.  

\subsection{Invariance} \label{s:invariance}

Let $J' \subset J$ denote the subideal generated by the ``$b$'' coefficients of the matrices in (\ref{e:b1})--(\ref{e:b4}).  To be precise, $J'$ is generated by:
\begin{itemize}
\item The elements $
\sum_{i=1}^r \boldsymbol{\epsilon}_i \boldsymbol{b}_i 
$
 for each row of type I in $D$.
\item The elements 
$
(\boldsymbol{a}_i + \boldsymbol{\nu_i})\boldsymbol{b}_j 
+ \boldsymbol{b}_i \boldsymbol{d}_j - \sum_{k = 1}^r \boldsymbol{\delta}_{ijk} \boldsymbol{b}_k
$ for each row of type II in $D$.
\item For each $v \in S \setminus \{v_0\}$ and distinct $\sigma, \tau \in \cB_v$, the elements
\[ B(\sigma, \tau) = \boldsymbol{b}_\sigma(\boldsymbol{x}_{\tau} - \boldsymbol{a}_\tau) - 
\boldsymbol{b}_\tau(\boldsymbol{x}_{\sigma} - \boldsymbol{a}_\sigma). \]
\end{itemize}

\begin{lemma} \label{l:stable}
The ideals $J'$ and $J$ are stable under the action of $\B$.  More precisely,  the $\Z$-module spanned by each set of 4 relations in (\ref{e:b1})--(\ref{e:b4}) is isomorphic as a $\B$-module to a copy of the adjoint representation $\cA$ (see Example~\ref{e:adjoint}).
\end{lemma}

\begin{proof}
This is clear for the relations in (\ref{e:b1}) and (\ref{e:b2}) since these relations are linear combinations of products of the $\boldsymbol{\rho}_i$, with coefficients in $R_0$ (on which $\B$ acts trivially), and the definition of our action is by simultaneous conjugation on  $\boldsymbol{\rho}_i$. We must verify this for (\ref{e:b4}) by direct computation.  The Borel $\B$ is generated by its torus and unipotent subgroups:
\begin{equation} \label{e:tu}
 \left\{ \sigma_{x,y} = \mat{x}{0}{0}{y} \right\}, \qquad \left\{ \tau_x = \mat{1}{0}{x}{1} \right\}. \end{equation}

Acting on the adjoint $\cA$, the element $\sigma_{x,y}$ fixes $A$ and $D$, scales $B$ by $y/x$ and scales $C$ by $x/y$.  From the definitions in (\ref{e:ABCD}), it follows that $\sigma_{x,y}$ acts the same way on $A(\sigma, \tau), \dots, D(\sigma, \tau)$.  Similarly, $\tau_x$ has the following action:
\[ \mat{A}{B}{C}{D} \mapsto \mat{A + Bx}{B}{C + (D-A)x - Bx^2 }{D-Bx},
\]
and one checks that the action on the matrix (\ref{e:b4}) is the same.  For instance,
\[
\tau_x B(\sigma, \tau)  = \boldsymbol{b}_\sigma(\boldsymbol{x}_{\tau} - \boldsymbol{a}_\tau- \boldsymbol{b}_\tau x) - 
\boldsymbol{b}_\tau(\boldsymbol{x}_{\sigma} - \boldsymbol{a}_\sigma - \boldsymbol{b}_\sigma x)
 = B(\sigma, \tau).
\]
The verification of the other 3 coefficients is similar and left to the reader.
The stability of $J'$ under $\B$ follows since the $b$-coefficient of the adjoint is stable under $\B$.
\end{proof}

The goal of the rest of this subsection is to prove the following.

\begin{lemma} We have $\overline{\boldsymbol{e}} \in H^0(\B, R/J')$. \label{l:ebar}
\end{lemma}

 The matrix $\boldsymbol{E}$ has coefficients in $R_0$, on which $\B$ acts trivially, so we must show that 
\[ \overline{\det(\boldsymbol{E}')} \in H^0(\B, R/J'). \]
Again we use the fact that $\B$ is generated by the $\sigma_{x,y}$ and $\tau_x$ defined in (\ref{e:tu}).  The matrix $\boldsymbol{E}'$ has coefficients in $R_0[\boldsymbol{a}_i, \boldsymbol{d}_i]$ on which  $\sigma_{x,y}$ acts trivially.  Therefore we must only consider the action of $\tau_x$.

Note that the first $t$ rows of $\boldsymbol{E}'$, the rows of type IV, and the rows of type V all have the same shape.  To each such row is attached a $v \in S\setminus\{v_0\}$ and a $\sigma \in \cG_v$.  
In the column corresponding to $v$ we have the element $\boldsymbol{x}_{\sigma} - \boldsymbol{a}_{\sigma}$ and in the other columns we have 0.  Therefore to streamline the exposition we will treat these rows the same and call them ``of type L.''  So $\boldsymbol{E}'$ has rows of type I, II, III, or L.

Conjugation by $\tau_x$ sends $\boldsymbol{a}_i \mapsto \boldsymbol{a}_i + \boldsymbol{b}_i x$ and $\boldsymbol{d}_i \mapsto \boldsymbol{d}_i - \boldsymbol{b}_i x$.  This action fixes rows of type I or III, but is non-trivial on rows of type II or L.  We will handle the rows of type II first, and then the rows of type L one $v$ at a time.  
Write $S \setminus\{v_0\} = \{v_1, \dots, v_{s+t} \}$.
For $i = 0, \dots, s+t$, let 
$\boldsymbol{E}'_i$  be the matrix obtained from $\boldsymbol{E}'$, 
but with the rows of type II and rows of type L associated to $v_j$ for $j \le i$ replaced by their image under $\tau_x$.
In particular $\boldsymbol{E}'_{s+t} = \tau_x(\boldsymbol{E}').$

\begin{lemma} We have $\det(\boldsymbol{E}'_0) \equiv \det(\boldsymbol{E}') \pmod{J'}$. \label{l:e0}
\end{lemma}

\begin{proof}
Suppose we have a row $w$ of type II associated to a pair $(i,j)$.  Then 
\[ \tau_x(w) - w = (0, 0, \dots, \boldsymbol{b}_j x, 0, \dots, -\boldsymbol{b}_i x, 0, \dots, 0), \]
where there is $\boldsymbol{b}_j x$ in the $\boldsymbol{\rho}_i$ slot and $ -\boldsymbol{b}_i x$ in the 
$\boldsymbol{\rho}_j$ slot.
The difference $\det(\boldsymbol{E}'_0) - \det(\boldsymbol{E}')$ is a linear combination (with coefficients $\pm 1$) 
of determinants of all matrices $M$ obtained by starting with $\boldsymbol{E}'$ and replacing some nonempty subset of the rows $w$ of type II with $\tau_x(w) - w$.

Suppose we are given such a matrix $M$ and row $w$.  We want to show that $\det(M) \equiv 0 \pmod{J'}.$  
Note that if $i = j$ then $\tau_x(w) - w = 0$ and hence $\det(M)=0$.  So we assume $i \neq j$.
We make the following alterations to $M$ which do not change the determinant:
\begin{itemize}
\item We replace $\tau_x(w) - w$ by 
\[ (0, 0, \dots,  x, 0, \dots, -x, 0, \dots, 0). \]
At the same time, in every row other than $w$ we multiply the $\boldsymbol{\rho}_j$-coordinate by $\boldsymbol{b}_j$ and the $\boldsymbol{\rho}_i$-coordinate by $\boldsymbol{b}_i$. 
\item We then add the new  $\boldsymbol{\rho}_j$ column to the new $\boldsymbol{\rho}_i$-column.
\item For each $1 \le k \le r, k \neq i, j$, we add $\boldsymbol{b}_k$ times the $\boldsymbol{\rho}_k$-column to the 
new $\boldsymbol{\rho}_i$-column.
\end{itemize}
Let us make some observations about the matrix $M$ that results from these changes:
\begin{itemize}
\item In rows of type I, II, or III, the $\boldsymbol{\rho}_i$-coordinate is precisely the generator of $J'$ associated to the row.
\item In rows of type L, we have $\boldsymbol{b}_\sigma$ in the $\boldsymbol{\rho}_i$-coordinate and 
$\boldsymbol{x}_\sigma - \boldsymbol{a}_\sigma$ in the $v$-coordinate.  The $2 \times 2$ determinant of any two such rows (associated to the same $v$) is one of the generators of $J'$.
\end{itemize}
As we now explain, the fact that $\det(M) \in J'$ follows from these two facts.  We compute $\det(M)$  as the sum
\begin{equation} \label{e:detm}
\sum_{\pi} \sgn(\pi) \prod_{w} M_{w, \pi(w)} \end{equation}
as $\pi$ ranges over all bijections between the set of rows and columns of $M$, and $w$ runs over the set of rows of $M$.
For a given $\pi$, consider the $w$ such that $\pi(w)$ is the  $\boldsymbol{\rho}_i$-column.  If $w$ is of type I, II, or III, then the contribution to $\det(M)$ lies in $J'$ by the first bulleted observation.  Now suppose  $w$ is of type L, associated to some $v \in S \setminus\{v_0\}$ and $\sigma \in \cG_v$.
The only rows of $M$ with nonzero component in the $v$-column are those of type L associated to $v$.  So we only obtain a nonzero contribution to $\det(M)$ if there is such a row $w'$, say attached to $\tau \in \cG_v$, and 
 $\pi(w') = $ the $v$-column.  In this case we consider the permutation $\pi'$ obtained from $\pi$ by swapping the roles of $w$ and $w'$.  So $\pi'(w) = v$-column, $\pi'(w') = \boldsymbol{\rho}_i$-column, and $\pi'(u) = \pi(u)$ for $u \neq w, w'$.  Then $\sgn(\pi') = - \sign(\pi)$ and by the second bulleted observation we find
 \[ \prod_{u} M_{u, \pi(u)} - \prod_{u} M_{u, \pi'(u)} = (\boldsymbol{b}_\sigma(\boldsymbol{x}_\tau - \boldsymbol{a}_\tau) - 
 \boldsymbol{b}_\tau(\boldsymbol{x}_\sigma - \boldsymbol{a}_\sigma))
  \prod_{u \neq w, w'} M_{u, \pi'(u)} \in J'. \]
 
Pairing off the permutations $(\pi, \pi')$ in this way shows that $\det(M) \in J'$. This concludes the proof.
\end{proof}

By a similar argument, we can show:

\begin{lemma}  For $i = 1, \dots, s+t$, we have $\det(\boldsymbol{E}'_i) \equiv \det(\boldsymbol{E}'_{i-1}) \pmod{J'}$.  \label{l:ei}
\end{lemma}

\begin{proof} Note that $\boldsymbol{E}'_i$ and $\boldsymbol{E}'_{i-1}$ only differ in the rows of type L associated to $v_i$, and only in the column associated to $v_i$.  In these entries, $\boldsymbol{E}'_{i-1}$ has the value $\boldsymbol{x}_\tau - \boldsymbol{a}_\tau$ and $\boldsymbol{E}'_i$
has the value $\boldsymbol{x}_\tau - \boldsymbol{a}_\tau - \boldsymbol{b}_\tau x$.  

We therefore modify $\boldsymbol{E}'_i$ by adding, for each $j=1, \dotsc, r$, the $\boldsymbol{\rho}_j$-column scaled by $\boldsymbol{b}_j x$ to the $v_i$-column.  The resulting matrix $M$ satisfies $\det(M) = \det(\boldsymbol{E}'_i)$.  Furthermore
the entries of $M$ and $\boldsymbol{E}'_{i-1}$ are equal except in the column associated to $v_i$ for the rows {\em not} associated to $v_i$.  In the $v_i$-column we have:
\begin{itemize}
\item For rows of type I, II, or III, precisely the associated generator of $J'$.
\item For rows of type L associated to $v_j \neq v_i$, the element $\boldsymbol{b}_\sigma x.$
\end{itemize}
Again we compute $\det(M)$ and $\det(\boldsymbol{E}'_{i-1})$ via (\ref{e:detm}), this time focusing on the row $w$ such that $\pi(w) = v_i$-column.  

If $w$ is of type I, II, or III, then the first bulleted point shows that the contribution to $\det(M)$ is 0 modulo $J'$. The same is true for 
$\det(\boldsymbol{E}'_{i-1})$ since the $(w, \pi(w))$-entry of $\boldsymbol{E}'_{i-1}$ is 0.

  If $w$ is of type L associated to $v_i$, then the $(u, \pi(u))$-entries of $M$ and $\boldsymbol{E}'_{i-1}$ are equal for all rows $u$, so the contributions to the determinants are the same. 
  
   Finally, if $w$ is of type L associated to some other place $v_j \neq v_i$, then the $(w, \pi(w))$-entry of $\boldsymbol{E}'_{i-1}$ is 0, so we must show that the contribution to $\det(M)$ is also $0$ modulo $J'$.  By the second bullet point above, the $(w, v_i)$-entry of $M$ is $\boldsymbol{b}_\sigma x$ and the $(w, v_j)$-entry is either $\boldsymbol{x}_\sigma - \boldsymbol{a}_\sigma$ or $\boldsymbol{x}_\sigma - \boldsymbol{a}_\sigma -
 \boldsymbol{b}_\sigma x$, depending on whether $j > i$ or $j \le i-1$.   Now, the only other rows of $M$ that have nonzero $v_j$-coordinate are those of type L associated to $v_j$. So the contribution of $\pi$ to $\det(M)$ is 0 unless $\pi(w') = v_j$-column for some other $w'$ of type L associated to $v_j$, say attached to the element $\tau \in \cG_{v_j}$.  The entries of $M$ in the $w'$ row are the same as those described for the $w$ row, with $\sigma$ replaced by $\tau$.  In both the cases $j > i$ or $j \le i-1$, the $2 \times 2$ determinant given by the $(w, w')$ rows and $(v_i, v_j)$ columns is equal to 
\begin{equation} \label{e:baj}
 \boldsymbol{b}_\sigma(\boldsymbol{x}_\tau - \boldsymbol{a}_\tau) - 
 \boldsymbol{b}_\tau(\boldsymbol{x}_\sigma - \boldsymbol{a}_\sigma) \in J' .
 \end{equation}
 Once again we define $\pi'$ to be the permutation  $\pi$ with the roles of $w, w'$ swapped.
As before, (\ref{e:baj}) implies that the sum of the contributions of $\pi$ and $\pi'$ to $\det(M)$ is 0 modulo $J'$.  Of course, the contribution of $\pi'$ to $\det(\boldsymbol{E}'_{i-1})$ is also 0.
Pairing off $(\pi, \pi')$ in this way concludes the proof.
\end{proof}

Lemmas~\ref{l:e0} and~\ref{l:ei} imply the desired result
\[ \det(\tau_x(\boldsymbol{E}')) = \det(\boldsymbol{E}'_{s+t}) \equiv \det(\boldsymbol{E}') \pmod{J'}. \]
This concludes the proof of Lemma~\ref{l:ebar}.

\subsection{A cascade of cohomology classes}
 
 In this section, we will assume the existence of a certain morphism of complexes (see Theorem~\ref{t:comm} below) and use it to prove the following.
 
 \begin{theorem} \label{t:push}
 Let $\iota \colon J' \rightarrow J$ be the inclusion, let $j \ge 1$, and let  
\[ \iota_*\colon H^j(\B, J') \longrightarrow H^j(\B, J) \]
be the induced map on rational cohomology groups.  Then $\iota_* = 0$.
 \end{theorem}

 \begin{theorem}\label{t:comm}
  There exists a morphism of complexes of $\B$-modules:
\begin{equation} \label{e:kozsul2} \begin{tikzcd}
0 \ar[r]  & C^r  \ar[r,"f_r"] \ar[d,"\iota_r"]
 & C^{r-1}   \ar[r, "f_{r-1}"] \ar[d,"\iota_{r-1}"]
 &  \cdots \ar[r,"f_2"] & C^1 \ar[r,"f_1"]  \ar[d,"\iota_1"] &  C^0 = J' \ar[d, "\iota_0",shift right=5]  \ar[r]
 & 0
 \\
 0 \ar[r] & D^r \ar[r,"g_r"] 
 & D^{r-1}   \ar[r, "g_{r-1}"]
 &  \cdots \ar[r,"g_2"] &  D^1 \ar[r,"g_1"] &  D^0 = J  \ar[r]
 & 0
 \end{tikzcd}
 \end{equation}
 such that the complex $C^{\bullet}$ is exact and every module $D^i$, for $i > 0$, is an acyclic $\B$-module. 
 \end{theorem}
 
 We assume Theorem \ref{t:comm} for now and use it to prove Theorem~\ref{t:push}.
 
 \begin{proof}[Proof of Theorem~\ref{t:push}]
The image of the map $f_1$ is $J'$.   Let \[ \alpha_j \in H^j(\B, J') = H^j(\B, \im(f_1)). \]
Our goal is to show $\iota_{0, *} \alpha_j = 0$.  Applying the coboundary in 
the long exact sequence in cohomology associated to 
\[ 0 \lra \ker(f_1) \lra C^1 \lra \im(f_1) \lra 0 \] to $\alpha_j$ 
 yields a class $\alpha_{j+1} \in H^{j+1}(\B, \ker(f_1))$. The class $\alpha_{j+1}$ represents the obstruction to lifting $\alpha_j$ to a class in $H^j(\B, C^1)$.  Writing $\ker(f_1) = \im(f_2)$, we can view \[ \alpha_{j+1} \in H^{j+1}(\B, \im(f_2)). \] Repeat the process above, using the coboundary in the long exact sequence associated to 
\[ 0 \lra \ker(f_2) \lra C^2 \lra \im(f_2) \lra 0 \]
to obtain $\alpha_{j+2} \in H^{j+2}(\B, \ker(f_2))$.  Continuing in this way we obtain 
\[ \alpha_{j+i} \in H^{j+i}(\B, \ker(f_{i})) = H^{j+i}(\B, \im(f_{i+1})) \] for $i = 0, \dots, r$.  Note $\alpha_{j+r} = 0$ since $\ker(f_r)=0$.

 For each $i = 0, \dotsc, r$, we define
\[ \beta_{j+i} = \iota_{i, *} \alpha_{j+i} \in H^{j+i}(\B, \im(g_{i+1})). \] We aim to prove that $\beta_j = \iota_{0, *} \alpha_j $ vanishes.  The obstruction to $\beta_{j+i} \in H^{j+i}(\B, \im(g_{i+1}))$  lifting to a class in $H^{j+i}(\B, D^{i+1})$ is precisely the image of $\beta_{j+i+1}$ 
in $H^{j+i+1}(\B, \ker(g_{i+1}))$.  Now, $\beta_{j+r} = 0$ since $\alpha_{j+r}=0$, and hence we conclude that $\beta_{j+r-1}$ lifts to a class in $H^{j+r-1}(\B, D^r)$.  However, $D^r$ is $\B$-acyclic, so $H^{j+r-1}(\B, D^r) =0$ and hence $\beta_{j+r-1}=0$.  Therefore, $\beta_{j+r-2}$ lifts to a class in  $H^{j+r-2}(\B, D^{r-1})$; again this cohomology group vanishes so $\beta_{j+r-2} = 0$.  This downward cascading continues and we obtain $\beta_i = 0 $ for all $i= j, \dotsc, j+r$.  In particular $\beta_j = 0$ as desired.
\end{proof}

In the remainder of the paper we will prove Theorem~\ref{t:comm}.  This will conclude the proof of Theorem~\ref{t:rll}.

\section{Construction of Complexes}

\subsection{Generators of $J'$} \label{s:generators}

Recall the ring $R$ defined in (\ref{e:r0def})--(\ref{e:rdef}) and the ideal $J' \subset R$ defined in \S\ref{s:invariance}.
We distinguish three types of generators of $J'$:
\begin{enumerate}
\item[(I)] The elements $
\sum_{i=1}^r \boldsymbol{\epsilon}_i \boldsymbol{b}_i 
$
 for each row of type I in $D$.
\item[(II)] The elements 
$
(\boldsymbol{a}_i + \boldsymbol{\nu}_i)\boldsymbol{b}_j 
+ \boldsymbol{b}_i \boldsymbol{d}_j - \sum_{k = 1}^r \boldsymbol{\delta}_{ijk} \boldsymbol{b}_k
$ for each row of type II in $D$.
\item[(III)] For each $v \in S \setminus \{v_0\}$ and distinct $\sigma, \tau \in \cB_v$, the elements
\[ B(\sigma, \tau) = \boldsymbol{b}_\sigma\boldsymbol{b}'_{\tau} - 
\boldsymbol{b}_\tau\boldsymbol{b}'_{\sigma}, \]
where \begin{equation} \label{e:bpdef}
\boldsymbol{b}'_{\sigma} := \boldsymbol{x}_{\sigma} - \boldsymbol{a}_{\sigma}. \end{equation}
\end{enumerate}

We partition $\{1,\ldots,r \}$ according to whether $\boldsymbol{b_i}$ is associated to a place $v \in S$. We find (after re-indexing) that we have $\boldsymbol{b}_1,\ldots,\boldsymbol{b}_k$ not associated to any place, followed by, for each $v \in S \setminus \{v_0\}$, $\{ \boldsymbol{b}_{\sigma} \colon \sigma \in \mathcal{B}_v \}$. Note that since $\boldsymbol{b}_{\sigma} = 0$ in $R$ for $\sigma \in \mathcal{B}_{v_0}$, we have non-zero variables $\boldsymbol{b}_1, \ldots, \boldsymbol{b}_k, \boldsymbol{b}_{k+1},\ldots,\boldsymbol{b}_{\ell}$ for some $k \le \ell \le r$.

We may view the generators of types I and II as linear forms in the variables $\boldsymbol{b}_1,\ldots, \boldsymbol{b}_{\ell}$: \begin{equation} L_i = \sum_{j=1}^{\ell} \boldsymbol{V}_{ij} \boldsymbol{b}_j \in R. \end{equation} The coefficients $\boldsymbol{V}_{ij} \in R$ are \emph{generic}, in that we can present $R$ as 
\[ R = S_0[\boldsymbol{a}_i, \boldsymbol{c}_i, \boldsymbol{d}_i]_{i=1,\ldots,r}[\boldsymbol{V}_{ij}][\boldsymbol{b}_1,\ldots, \boldsymbol{b}_{\ell}], \] 
where $S_0$ is a subring of $ R_0 = \Z[ \boldsymbol{\nu}_i, \boldsymbol{\epsilon}_{i}, \boldsymbol{\delta}_{ijk}, \boldsymbol{x}_\sigma]$ and the $\boldsymbol{V}_{ij}$ are certain linear combinations of the generators of $R_0$.
Specifically, each $\boldsymbol{V}_{ij}$ is an expression of the form $\boldsymbol{\epsilon}_i $, $(\boldsymbol{a}_i + \boldsymbol{\nu}_i) - \boldsymbol{\delta}_{ijj}, \boldsymbol{d}_j - \boldsymbol{\delta}_{iji},$ or $-\boldsymbol{\delta}_{ijk}$, and $S_0$ is obtained by deleting (in each of the 4 respective cases) $\boldsymbol{\epsilon}_i, \boldsymbol{\delta}_{ijj}, \boldsymbol{\delta}_{iji},$ or $\boldsymbol{\delta}_{ijk}$ from the definition of $R_0$.
 Note that the group $\B$ acts trivially on the $\boldsymbol{V}_{ij}$'s of the form $\boldsymbol{\epsilon}_i $ or $-\boldsymbol{\delta}_{ijk}$, but acts non-trivially on $(\boldsymbol{a}_i + \boldsymbol{\nu}_i) - \boldsymbol{\delta}_{ijj}$ and  $\boldsymbol{d}_j - \boldsymbol{\delta}_{iji}$. Thus for the purpose of producing the resolutions of Theorem \ref{t:comm}, we may consider these generators simply as $\sum \boldsymbol{V}_{ij} \boldsymbol{b}_j$, but for the purpose of proving $\B$-acyclicity we will need to distinguish these two types of generators.

We may do a further change of variables and present $R$ as 
\begin{equation}\label{e:present-R}
R = T_0[\boldsymbol{b}'_{\sigma} , \boldsymbol{V}_{ij}][\boldsymbol{b}_1,\ldots,\boldsymbol{b}_{\ell}],
\end{equation}
for $T_0$ a subring of $R_0[\boldsymbol{a}_i, \boldsymbol{c}_i, \boldsymbol{d}_i]_{i=1,\ldots,r}$.  Here the $\boldsymbol{b}'_{\sigma}$ are defined in (\ref{e:bpdef}).

\subsection{Motivation via regular sequences}

For each $v \in S \setminus \{ v_0 \}$, fix $\sigma_v \in \mathcal{B}_v$. Over the ring $R' = R[(\boldsymbol{b}'_{\sigma_v})^{-1}]_{v \in S \setminus \{ v_0 \}}$, we may choose a smaller set of elements to generate the extended ideal $R' J'$, by removing certain generators of type III. Indeed, the following elements generate $R'J'$:
\begin{itemize}
\item the linear forms $\sum \boldsymbol{V}_{ij} \boldsymbol{b}_j$,
\item for each $v \in S \setminus \{ v_0 \}$, $B(\sigma_v,\tau) =  \boldsymbol{b}'_{\tau}\boldsymbol{b}_{\sigma_v} - 
\boldsymbol{b}'_{\sigma_v} \boldsymbol{b}_\tau$ for all $\tau \in \cB_v \setminus \{ \sigma_v \}$.
\end{itemize}

Indeed, for $\sigma \in \mathcal{B}_v$, we have 
\[ B(\sigma, \tau) = (\boldsymbol{b}'_{\sigma_v})^{-1}( \boldsymbol{b}'_{\sigma} B(\sigma_v, \tau) - \boldsymbol{b}'_{\tau} B(\sigma_v, \sigma)). \]
These bulleted elements form a regular sequence in $R'$, a fact whose whose proof we omit as it will not end up being used. To obtain an exact complex of $R'$-modules resolving $R'/R' J'$ as required in Theorem~\ref{t:comm}, we would simply take the Koszul complex on this regular sequence.

However, to prove Theorem \ref{t:comm} we need to construct a resolution of $R/J'$, not of $R'/R'J'$. Unfortunately, the ring $R/J$ cannot typically be resolved using a Koszul complex, as the ideal $J'$ is not typically generated by a regular sequence. Nevertheless, it is still possible to resolve $R/J'$ using a generalization of the Koszul complex due to Buchsbaum--Rim \cite{b}, which we now describe.

\subsection{Buchsbaum--Rim complexes}
For this section, let $R$ denote any commutative ring. Associated to an $R$-linear map \[ f \colon V \lra W,  \qquad V = R^n,\ W = R^m, m \leq n, \] Buchsbaum--Rim \cite{b} define two complexes $R(f)$ and $R(\det(f))$.

\subsubsection{Special cases}

\noindent \textbf{Case $m = 1$.}
For a map $f \colon V \lra R$, \[ R(f) = R(\det(f)) = \text{Koszul complex on the elements } f(e_1),\ldots, f(e_n) \in R. \] The degree $k$ term of the complex $R(f)$ is $\bigwedge^k V$, and in degrees $1$ and $0$ the complex is given by $\begin{tikzcd} V \ar[r,"f"] & R \end{tikzcd}$.  Here and throughout, all alternating powers are understood to be over $R$.

\bigskip

\noindent\textbf{Case $m = 2$.}
Consider a map $f \colon V \lra W$, with \[ V = R^n = \bigoplus_i R e_i, \qquad W = R^2 = R w_1 \oplus R w_2. \] Write $f(e_i) = b_i w_1 + b_i' w_2$, and define $r_{ij} := b_i b'_j - b_j b'_i$.
The complex $R(f)$ begins as follows:
\begin{equation}
\begin{tikzcd}
 \cdots \ar[r] & \bigwedge^2(W^*) \otimes \bigwedge^3 V \ar[r] &  V \ar[r, "f"] & W.
\end{tikzcd}
\end{equation}
The image of the map $\bigwedge^2(W^*) \otimes \bigwedge^3 V \lra V$ is generated as an $R$-module by  \[ d_{ijk} = r_{ij} e_k + r_{jk} e_i + r_{ki} e_j. \]
The complex $R(\det(f))$ begins as follows:
\begin{equation}
\begin{tikzcd}
 \cdots \ar[r] & \Big(W^* \otimes \bigwedge^3 V\Big) \oplus \Big(\bigwedge^2(W^*) \otimes \bigwedge^4 V\Big) \ar[r] & \bigwedge^2 V \ar[r, "\det(f)"] & \bigwedge^2 W.
\end{tikzcd}
\end{equation}

For $k \geq 1$, the degree $k$ term of the complex $R(\det(f))$ is:
\begin{equation} \bigoplus_{s_i \in \{1, 2\}} \Big(\bigwedge\nolimits^{s_1} W^*\Big) \otimes \cdots \otimes \Big(\bigwedge\nolimits^{s_{k-1}} W^*\Big) \otimes \bigwedge\nolimits^{2 + \sum s_i} V. \end{equation}

The image of the map $\det(f) \colon \bigwedge^2(W^*) \otimes \bigwedge^2 V \lra R$ is generated over $R$ by the $2 \times 2$ minors of the $m \times 2$ matrix corresponding to $f$.

The complex $R(\det(f))$ is functorial in the following sense. Given another $R$-linear map $f' \colon V' \lra W$, $V' = R^{n'}$, as well as an $R$-linear map $g \colon V \lra V'$ such that $f \circ g = f'$, there is a chain homomorphism $\det(g) \colon R(\det(f)) \lra R(\det(f'))$ defined using the map $\bigwedge^i g \colon \bigwedge^i V \lra \bigwedge^i V'$.  

\subsubsection{Exactness of BR-complexes}
For $m \leq k \leq n$, define $f(k) \colon R^k \lra R^m$ by restricting $f$ to the first $k$ summands of $R^n$. 
\begin{definition}[\cite{b}]
The map $f$ is {\em regular} if for all $m \leq k \leq n$, we have $H_1(R(f(k))) = 0$. 
\end{definition}

The following result of Buchsbaum--Rim generalizes the well-known exactness of the Kozsul complex associated to a regular sequence.
\begin{theorem}[\cite{b}]  The complexes $R(f)$ and $R(\det(f))$ satisfy the following.
\begin{enumerate}
\item If $f$ is regular then $H_i(R(f)) = 0$ for all $i \geq 1$.
\item If $H_i(R(f)) = 0$ for $i \geq j$, then $H_i(R(\det(f))) = 0$ for $i \geq j$.
\end{enumerate}
\end{theorem}

\begin{corollary}
If $f \colon R^n \lra R^m$ is regular then $R(\det(f))$ is a resolution of $R/\im(\det(f))$.
\end{corollary}

We will apply this result to multiple complexes $R(\det(f))$ and then take the tensor product, necessitating the following lemma.

\begin{lemma}\label{l:tensor}
Consider a ring $R$ with maps $f_i \colon R^{n_i} \lra R^{m_i}$, $1\leq i \leq k$, such that $m_i \leq n_i$. Define ideals $J_i = \im(\det(f_i))$, and let $J = J_1 + \cdots + J_{k}$.
 Suppose that $f_i \otimes R/(J_1 + \cdots + J_{i-1})$ is regular for all $i \leq k$. Then the tensor product over $R$ of these complexes, $\bigotimes_{j=1}^k R(\det(f_j))$, is a resolution of $R/J$.
\end{lemma}
\begin{proof}
By induction on $k$, we may suppose that the natural map \[ \otimes_{j=1}^{k-1} R(\det(f_j)) \lra R/\Big(\sum_{j=1}^{k-1} J_j\Big)\] is a quasi-isomorphism. Then, since $R(\det(f_k))$ is a complex of free $R$-modules, the Tor spectral sequence (see \cite{stacks-tor}) implies that 
\[ \otimes_{j=1}^{k} R(\det(f_j)) \lra R/\Big(\sum_{j=1}^{k-1} J_j\Big) \otimes R(\det(f_k)) \] is also a quasi-isomorphism. By the assumption on regularity, the codomain is an exact complex such that the natural map \[ R/\Big(\sum_{j=1}^{k-1} J_j\Big) \otimes R(\det(f_k)) \lra R\Big(\sum_{j=1}^{k} J_j\Big) \] is a quasi-isomorphism. 
\end{proof}

\subsection{Outline of proof of Theorem \ref{t:comm}} \label{s:outline}

Recall that in \S\ref{s:generators} we described generators for $J'$ of the form:
\begin{itemize}
\item $L_i = \sum \boldsymbol{V}_{ij} \boldsymbol{b}_j$, $i = 1, \dotsc, n$.
\item $B(\sigma, \tau) = \boldsymbol{b}_\sigma\boldsymbol{b}'_{\tau} - 
\boldsymbol{b}_\tau\boldsymbol{b}'_{\sigma}$ for $v \in S \setminus \{v_0\},$ distinct $\sigma, \tau \in \cB_{v_0}$.
\end{itemize}
Recall also the $\mathbf{B}$-modules $\cB$ and $\cV$ defined in Example \ref{e:adjoint}, and the notation
   $R(1) = \cB \otimes_{\Z} R.$ We also write $\cV_R = \cV \otimes_{\Z} R$; this is $\cong R^2$ as an $R$-module.

To handle the first collection of generators bulleted above, we define the $\B$-equivariant map
\begin{equation} \label{e:fdef}
 f \colon R(1)^n = \bigoplus_{i=1}^n R(1)e_i \lra R, \qquad f(e_i) = L_i.  
 \end{equation}
For the second collection of generators bulleted above, we define for each $v \in S \setminus \{v_0\}$ the 
$\B$-equivariant map 
\begin{equation} f_v \colon \bigoplus_{\sigma \in \cB_{v}} R(1)e_{\sigma} \lra \mathcal{V}_R, \qquad
f_v(e_{\sigma}) = A \otimes \boldsymbol{b}_{\sigma} + B \otimes\boldsymbol{b}'_{\sigma}.
\end{equation}

We consider the Buchsbaum--Rim complexes associated to the maps $f, f_v$:
\begin{itemize}
\item The Koszul complex $R(f)$ is given in degrees $1$ and $0$ by the map $f$, whose image is the ideal $J'_0 := (\{\sum \boldsymbol{V}_{ij} \boldsymbol{b}_j\}_i)$.
\item The Buchsbaum--Rim complex $R(\det(f_v))(-1)$ (note the twist) is given in degrees $1$ and $0$ by the map \[ \det(f_v) \colon \Big(\bigwedge^2 \oplus_{\sigma} R(1)\Big) \otimes R(-1) \lra \bigwedge^2\mathcal{V}_R \otimes R(-1) \cong R. \] The image of $\det(f_v)$ is the ideal $J'_v := (B(\sigma,\tau) \colon \sigma, \tau \in \mathcal{B}_v) \subset R$. 
\end{itemize}

The tensor product of these complexes, \begin{equation} C^{\bullet} := R(f) \otimes_R \bigotimes_{v \in S \setminus \{ v_0 \}} R(\det(f_v))(-1), \label{e:cdef}
 \end{equation} is a complex of $\B$-modules such that $H^0(C^{\bullet}) = R/(J'_0 + \sum_v J'_v) = R/J'$.

To prove Theorem \ref{t:comm} it will suffice to prove:
\begin{prop}\label{p:resolve-j} The complex $C^{\bullet}$ satisfies the following.
\begin{enumerate}
\item $H^i(C^{\bullet}) = 0$ for $i > 0$, i.e. $C^{\bullet}$ is a resolution of $R/J'$ and the image of $C^1 \lra C^0 = R$ is the ideal $J'$.
\item $C^{\bullet}$ maps to a complex $D^{\bullet}$ such that the image of $D^1 \lra D^0 = R$ is $J$, and $D^i$ is an acyclic $\B$-module all $i \geq 0$.
\end{enumerate}
\end{prop}

We now explain how this proposition follows from the results of \S\ref{ss:extend-f}-\S\ref{ss:tensor-exact}. In the next two sections (\S \ref{ss:extend-f} and  \S \ref{ss:extend-fv}), we will extend the maps $f$ and $f_v$ to maps \[ \widetilde{f} \colon \bigoplus_{i=1}^n \mathcal{A}_R \lra R,\qquad \widetilde{f}_v \colon \bigoplus_{\sigma \in \mathcal{B}_v} \mathcal{V}_R \lra \mathcal{V}_R. \]
Here $\cA_R = \cA \otimes_\Z R$ where $\cA$ denotes the adjoint representation
defined in example \ref{e:adjoint}. 
There are natural maps $R(f) \lra R(\widetilde{f})$ and $R(\det(f_v))(-1) \lra R(\det(\widetilde{f}_v))(-1)$. We will show that these maps factor through subcomplexes \[ R(\widetilde{f})' \subset R(\widetilde{f}), \qquad R(\det(\widetilde{f}_v))(-1)' \subset R(\det(\widetilde{f}_v))(-1) \] consisting of good $\G$-modules tensored with the $\B$-module $R$. 

Define \[ J_0 := \text{ideal of } R  \text{ generated by the coefficients of (\ref{e:b1}) and (\ref{e:b2}),} \] and for $v \in S \setminus \{v_0\}$ define \[ J_v :=  (A(\sigma,\tau),\ B(\sigma,\tau),\ C(\sigma,\tau),\ D(\sigma,\tau) \colon \sigma, \tau \in \mathcal{B}_v) \subset R. \] Then by definition, \[ J = J_0 + \sum_{v \in S \setminus \{v_0\}} J_v. \]

We will show that $H^0(R(\widetilde{f})') = R/J_0$ and $H^0(R(\widetilde{f}_v)') = R/J_v$.
Define \begin{equation}D^{\bullet} := R(\widetilde{f})' \otimes_R \bigotimes_{v \in S \setminus \{ v_0 \}} R(\det(\widetilde{f}_v))(-1)'.  \label{e:ddef}
\end{equation} It follows that \[ H^0(D^{\bullet}) = R/(J_0 + \sum_{v \in S \setminus \{v_0\}} J_v) = R/J, \] i.e. that the image of the map $D^1 \lra D^0 = R$ is $J$. As all terms of $D^{\bullet}$ consist of good $\G$-modules tensored with $R$, we conclude by Theorem \ref{t:r-acyclic} that all terms of $D^{\bullet}$ are acyclic $\B$-modules.

It will then remain to show that $H^i(C^{\bullet}) = 0$ for $i > 0$. This will be proved using Proposition \ref{p:br-exact} below, building on the work of sections \S\ref{ss:br-regular}-\S\ref{ss:tensor-exact}.

\subsection{Extension of $f$}\label{ss:extend-f}

We define a map \[ \widetilde{f} \colon \cA_R^n =  \bigoplus_{i=1}^n \cA_R e_i \lra R \] extending the map $f$ defined in (\ref{e:fdef}) as follows (here {\em extension} refers to the natural inclusion $R(1) = \cB_R \subset \cA_R$, see Example~\ref{e:adjoint}).
\begin{itemize}
\item For linear forms $L_i$ corresponding to relations of type I, we set
 \begin{align*}
  \widetilde{f}(Ae_i) &= \sum \boldsymbol{V}_{ij} \boldsymbol{a}_j, \quad  \widetilde{f}(Be_i) = \sum \boldsymbol{V}_{ij} \boldsymbol{b}_j = L_i,\\
  \widetilde{f}(Ce_i) &= \sum \boldsymbol{V}_{ij} \boldsymbol{c}_j, \quad \widetilde{f}(D e_i) = \sum \boldsymbol{V}_{ij} \boldsymbol{d}_j. 
  \end{align*}
 In other words, as $L_i$ came from the $b$-coefficient of some matrix $\sum_{i=1}^r \boldsymbol{\epsilon}_i \boldsymbol{\rho}_i$, we extend $f$ using the other 3 coefficients of $\sum_{i=1}^r \boldsymbol{\epsilon}_i \boldsymbol{\rho}_i$.

\item For linear forms $L_i$ corresponding to relations of type II, i.e.\ to the $b$-coefficient of $(\boldsymbol{\rho}_i + \boldsymbol{\nu}_i) \boldsymbol{\rho}_j - \sum_{k = 1}^r \boldsymbol{\delta}_{ijk} \boldsymbol{\rho}_k$, we similarly extend $f$ using the other 3 coefficients.
 
\end{itemize}
 
The Koszul complex $R(\widetilde{f})$ is a complex of $\B$-modules. By the functoriality of the Koszul complex we obtain a $\B$-equivariant map $R(f) \lra R(\widetilde{f})$.

\begin{lemma} Expanding the terms of the complex $R(\widetilde{f})$ via the Kunneth formula, 
the image of $R(f) \lra R(\widetilde{f})$ is contained in the subcomplex $R(\widetilde{f})'$ whose $k$th term is 
\[\bigoplus_{\# I = k} \bigotimes_{i \in I} (\cA e_i) \otimes R. \]
\end{lemma}
\begin{proof}
The $k$th degree term of the Koszul complex $R(f)$ is given by \[\bigwedge^k \Big(\bigoplus_i R(1) e_i \Big) \cong \bigoplus_{\# I = k} \bigotimes_{i \in I} R(1), \] 
since $R(1)$ is free of rank 1 over $R$. The image of this map is contained in the subcomplex $R(\widetilde{f})'$ whose $k$th degree term is identified by the Kunneth formula with \[ \bigoplus_{\# I = k} \bigotimes_{i \in I} \cA_R e_i, \]
i.e.\ the submodule of $\bigwedge^k\bigoplus \cA_R e_i$ generated by wedge products $v_{i_1} e_{i_1} \wedge \cdots \wedge v_{i_k} e_{i_k}$, $v_{i_j} \in \mathcal{A}_R$, such that $i_1 < \cdots < i_k$.
\end{proof}

\begin{lemma}\label{l:adj-r-acyclic}
The terms of $R(\widetilde{f})'$ are the tensor product of good $\G$-modules with $R$.
\end{lemma}
\begin{proof}
This follows directly from Corollary \ref{c:adj-good}.
\end{proof}

\subsection{Extension of $f_v$}\label{ss:extend-fv}

For $\sigma \in \cB_v$, recall the notation $\boldsymbol{b'_\sigma} = \boldsymbol{x_\sigma} - \boldsymbol{a_\sigma}$. Define
now also $\boldsymbol{c}'_\sigma = \boldsymbol{x}_\sigma - \boldsymbol{d}_\sigma$.
We define \[ \widetilde{f}_v \colon \bigoplus_{\sigma \in \cB_v} \mathcal{V}_R e_{\sigma} \lra \mathcal{V}_R\] by \begin{align*}
\widetilde{f}_v(A e_{\sigma}) &= A \otimes \boldsymbol{c'_\sigma} + B \otimes \boldsymbol{c_\sigma}, \\
\widetilde{f}_v(B e_{\sigma}) &= A \otimes \boldsymbol{b_\sigma} + B \otimes \boldsymbol{b'_\sigma}. 
\end{align*}

\begin{lemma}
The map $\widetilde{f}_v$ is $\B$-equivariant.
\end{lemma}
\begin{proof}

Let $g = \mat{x}{0}{y}{z}$. We have
\begin{align*}
 g \cdot \boldsymbol{b_{\sigma}} &= \frac{z}{x} \boldsymbol{b_{\sigma}}, \\
 g \cdot \boldsymbol{b'_{\sigma}} &= g \cdot (\boldsymbol{x_{\sigma}} - \boldsymbol{a_{\sigma}}) = \boldsymbol{x_{\sigma}} - \boldsymbol{a_{\sigma}} - \frac{y}{x} \boldsymbol{b_{\sigma}} = \boldsymbol{b'_{\sigma}} - \frac{y}{x} \boldsymbol{b_{\sigma}} \\
 g \cdot \boldsymbol{c_{\sigma}} &= -\frac{y}{z}\boldsymbol{a_{\sigma}} - \frac{y^2}{xz}\boldsymbol{b_{\sigma}} + \frac{x}{z}\boldsymbol{c_{\sigma}} + \frac{y}{z} \boldsymbol{d_{\sigma}}\\
 g \cdot \boldsymbol{c'_{\sigma}} &= g \cdot (\boldsymbol{x_{\sigma}} - \boldsymbol{d_{\sigma}}) = \boldsymbol{x_{\sigma}} - \boldsymbol{d_{\sigma}} + \frac{y}{x} \boldsymbol{b_{\sigma}} =  \boldsymbol{c'_{\sigma}} + \frac{y}{x} \boldsymbol{b_{\sigma}}. 
 \end{align*}
We compute
\[ g \cdot \left(A \otimes \boldsymbol{b_{\sigma}} + B \otimes \boldsymbol{b'_{\sigma}}\right)  = \left(A + \frac{y}{x} B\right)  \otimes \left(\frac{z}{x} \boldsymbol{b_{\sigma}}\right)  + \left(\frac{z}{x} B\right)  \otimes  \left(\boldsymbol{b'_{\sigma}} - \frac{y}{x} \boldsymbol{b_{\sigma}}\right)  = \frac{z}{x}\left( A \otimes \boldsymbol{b_{\sigma}} + B \otimes \boldsymbol{b'_{\sigma}}\right) . \]

We also have
\begin{align*}
g \cdot\left(A \otimes \boldsymbol{c'_{\sigma}} + B \otimes \boldsymbol{c_{\sigma}} \right) &=  \left(A + \frac{y}{x} B\right)  \otimes \left(\boldsymbol{c'_{\sigma}} + \frac{y}{x} \boldsymbol{b_{\sigma}}\right)  + \left(\frac{z}{x} B\right)  \otimes \left(-\frac{y}{z}\boldsymbol{a_{\sigma}} - \frac{y^2}{xz}\boldsymbol{b_{\sigma}} + \frac{x}{z}\boldsymbol{c_{\sigma}} + \frac{y}{z} \boldsymbol{d_{\sigma}}\right)  \\
 &= A \otimes \left( \left(\boldsymbol{x_{\sigma}} - \boldsymbol{d_{\sigma}}\right)  + \frac{y}{x} \boldsymbol{b_{\sigma}}\right)  \\ & \ \ \ +  B\otimes \left( \frac{y}{x}\left(\boldsymbol{x_{\sigma}} - \boldsymbol{d_{\sigma}}\right)  + \frac{y^2}{x^2} \boldsymbol{b_{\sigma}} - \frac{z}{y} \boldsymbol{a_{\sigma}} - \frac{z}{x}\boldsymbol{b_{\sigma}} + \boldsymbol{c_{\sigma}}  + \frac{y}{x} \boldsymbol{d_{\sigma}}\right) \\
 &= \left(A \otimes \boldsymbol{c'_{\sigma}} + B \otimes \boldsymbol{c_{\sigma}} \right)  + \frac{y}{x}  \left( A \otimes \boldsymbol{b_{\sigma}} + B \otimes \boldsymbol{b'_{\sigma}}\right).
\end{align*}
The result follows.
\end{proof}

The BR-complex $R(\det(f_v))$ is a complex of $\B$-modules whose degree $0$ term is \[ \bigwedge\nolimits^2 \mathcal{V}_R \cong R(1),\] and whose degree $k \geq 1$ term is:
\begin{equation} \bigoplus_{s_i \in \{1, 2\}} \Big(\bigwedge\nolimits^{s_1} \mathcal{V}_R^*\Big) \otimes \cdots \otimes \Big(\bigwedge\nolimits^{s_{k-1}} \mathcal{V}_R^*\Big) \otimes \bigwedge\nolimits^{2 + \sum s_i} \Big(\bigoplus_{\sigma \in \cB_v} R(1)e_{\sigma}\Big). \end{equation}
The BR-complex $R(\det(\widetilde{f}_v))$ is the same but with $\bigwedge^{k} (\bigoplus_{\sigma} R(1)e_{\sigma})$ replaced by $\bigwedge^{k} (\bigoplus_{\sigma} \mathcal{V}_Re_{\sigma})$. The map $R(\det(f_v)) \lra R(\det(\widetilde{f}_v))$ has image contained in the subcomplex $R(\det(\widetilde{f}_v))'$ where each $\bigwedge^k (\bigoplus_{\sigma} \mathcal{V}_Re_{\sigma})$ is replaced by the $R$-submodule generated by wedge products $v_{\sigma_1} e_{\sigma_1} \wedge \cdots \wedge v_{\sigma_{k}} e_{\sigma_k}$, $v_i \in \mathcal{V}_R$, such that $\sigma_{i} \neq \sigma_{j}$ for $i \neq j$. This submodule is isomorphic to a direct sum of modules $\bigotimes^k_R \mathcal{V}_R$.

\begin{lemma}
$H^0(R(\det(\widetilde{f}_v))'(-1)) = R/J_v$.
\end{lemma}
\begin{proof}

We are restricting the map \[ \det(\widetilde{f}_v) \colon \bigwedge^2 (\oplus_{\sigma} \mathcal{V}_R e_{\sigma}) \otimes R(-1) \lra R \] to the wedge-products $v_{\sigma} e_{\sigma} \wedge v_{\tau} e_{\tau}$, $v_{\sigma}, v_{\tau} \in \mathcal{V}_R$, with $\sigma \neq \tau$. We need to show that the image of this restricted map is $J_v$.
A direct computation using the definitions in (\ref{e:ABCD}) shows that the elements $B e_{\sigma} \wedge A e_{\tau}$, $A e_{\sigma} \wedge A e_{\tau}$, $B e_{\sigma} \wedge B e_{\tau}$, $A e_{\sigma} \wedge B e_{\tau}$, map to $\pm A(\sigma, \tau)$, $\pm B(\sigma, \tau)$, $\pm C(\sigma, \tau)$, $\pm D(\sigma,\tau)$. By definition, these 
latter elements generate $J_v$.  The result follows.
\end{proof}

\begin{lemma} \label{l:localgood}
The terms of $R(\det(\widetilde{f}_v))'(-1)$ are the tensor products of good $\G$-modules with $R$.
\end{lemma}
\begin{proof}

The degree $0$ term of $R(\det(\widetilde{f}_v))'(-1)$ is $R$.  
The degree $k \geq 1$ term of $R(\det(\widetilde{f}_v))'(-1)$
 is a direct sum of modules of the form \begin{equation}\label{e:good1}\left(\Big(\bigwedge\nolimits^{s_1} \mathcal{V}_R^*\Big) \otimes \cdots \otimes \Big(\bigwedge\nolimits^{s_{k-1}} \mathcal{V}_R^*\Big) \otimes \mathcal{V}_R^{\otimes 2 + \sum s_i}\right) \otimes R(-1).\end{equation}

 Setting $m = \sum (s_i - 1) = \Big(\sum s_i\Big) - (k-1)$, (\ref{e:good1}) is isomorphic to 
 \begin{equation} \label{e:good2}
 \left(\Big(\bigwedge\nolimits^2 \mathcal{V}_R^* \Big)^{\otimes m} \otimes (\mathcal{V}_R^*)^{\otimes (k-1) - m} \otimes \mathcal{V}_R^{k+m+1}\right) \otimes R(-1)
 \end{equation}
 
 Using the isomorphisms $\bigwedge^2 \mathcal{V} \cong \Z(1)$, $\mathcal{V} \otimes \mathcal{V}^* \cong \mathcal{A}$, and $\mathcal{V} \otimes \mathcal{V} = \mathcal{A}(1)$, 
  one readily checks that (\ref{e:good2}) is isomorphic to $\mathcal{A}^{\otimes k} \otimes R$.
The result now follows from Corollary~\ref{c:adj-good}.
\end{proof}

As described in \S\ref{s:outline} above, Lemmas \ref{l:adj-r-acyclic} and \ref{l:localgood} imply that the complex $D^\bullet$
defined in (\ref{e:ddef}) satisfies the property that $D^i$ is $\B$-acyclic for $i > 0$.  In other words, we have finished the proof of the second part of Proposition~\ref{p:resolve-j}.  It remains to prove the first part of the proposition.

\subsection{Regularity criterion for maps $f \colon R^n \lra R^2$}\label{ss:br-regular}

Let $R$ be a commutative ring. Given a map $f \colon R^n \lra R^2$, write $f(e_i) = (b_i, b_i')$, and define $r_{ij} = b_i b'_j - b_j b'_i$.

\begin{lemma}\label{l:reg}
Suppose that for $k = 2, \dots, n$, the congruence
\[ r_{1k} x \equiv 0 \pmod{(r_{12}, \dots, r_{1(k-1)})} \]
implies that
\[ x \in (r_{ij} \colon i,j < k ). \]
Then the map $f \colon R^n \lra R^2$ is regular.
\end{lemma}
\begin{proof}
The goal is to show that the kernel of $f$ is generated as an $R$-module by \[ d_{ijk} = r_{ij}e_k + r_{jk}e_i + r_{ki}e_j. \] 

We proceed via induction. In the base case $n=2$, our assumption is that  $r_{12}$ is a non-zerodivisor, and this implies that $\ker(f) = 0$.  We proceed to the general case.

Any element $x = (x_1, \ldots, x_n) \in \ker(f)$ satisfies $\sum r_{1i} x_i = 0$. 
Hence \[ r_{1n}x_n \equiv 0 \mod (r_{12}, \ldots, r_{1(n-1)}). \] The assumption implies that  we may write $x_n = \sum_{i<j < n} \alpha_{ij} r_{ij}$ for some $\alpha_{ij} \in R$. Therefore the final coordinate of $(x_1, \ldots,x_n) - \sum \alpha_{ij} d_{ijn}$ is zero. Writing this vector as $y = (y_1, \dots, y_{n-1}, 0)$, we can apply the inductive hypothesis to $y$ and the restriction of $f$ to its first $n-1$ coordinates.  Hence $y$ is a linear combination of the $d_{ijk}$, so $x$ is as well.
\end{proof}

\begin{corollary} \label{c:genericb}
Let $R_0$ be a  commutative ring, and let $R = R_0[b'_1, \ldots, b'_n, b_1, \ldots, b_n]$. The map $f \colon R^n \lra R^2$ defined by $f(e_i) = (b_i, b_i')$ is regular.
\end{corollary}
\begin{proof}
We will show that the criterion of Lemma \ref{l:reg} holds. 
Let $J_k = (\{ r_{ij} \}_{1 \leq i,j \leq k}) \subset R$. 

First we write $V = b'_1$ and work in $R[V^{-1}]$.
Note that $R[V^{-1}]\cdot J_k = (r_{12},\ldots,r_{1k})$, and that \begin{equation} \label{e:rv}
 R[V^{-1}]/R[V^{-1}]\cdot J_k \cong R_0[b_1', \dots, b_n', b_1, b_{k+1}, \dots, b_n][(b_1')^{-1}]. \end{equation} 
 It is clear that the image of $r_{1(k+1)}$ is a non-zerodivisor in (\ref{e:rv}), and hence that the criterion of Lemma \ref{l:reg} holds.
 
It remains to prove that $R \cap (R[V^{-1}] \cdot J_k) = J_k.$  Therefore suppose that $x \in R$ and that the image of $x$ vanishes in $(\ref{e:rv})$.  Fix nonegative integers $c_1, \dots, c_k$, as well as nonnegative integers $a_i, a_i'$ for $k < i \le n$.
Applying the isomorphism (\ref{e:rv}), we see that the coefficients of the monomials in $x$ of the form
\[ b_1^{a_1} \cdots b_n^{a_n}(b_1')^{a_1'} \cdots (b_n')^{(a_n')} \]
with $a_i + a_i' = c_i$ for $i \le k$, $\sum a_i = \text{constant}$, and $\sum a_i' = \text{constant}$, must sum to zero.  For example, $x$ could be a linear combination of expressions of the form $b_2 b'_3 - b'_2 b_3$ or $d_1 b_2 b_3 b'_4 + d_2 b_2 b'_3 b_4 + d_3 b'_2 b_3 b_4$ with $d_1 + d_2 + d_3 = 0$. It is easy to see that, modulo $J_k$, all such monomials are equivalent.  Therefore $x \in J_k$.
\end{proof}

\subsection{Regularity criterion for sequences of linear polynomials}\label{ss:regular-seq}

In this subsection we will show that the regularity criterion above hold for generic linear polynomials.  As motivation, 
let $k$ be a field and suppose we are given homogeneous linear polynomials \[ L_i = \sum_{j=1}^n a_{ij} X_j \in k[X_1, \ldots, X_n], \quad 1 \leq i \leq m. \] If the $m \times n$ matrix $A = (a_{ij})$ has linearly independent rows, then there is a change of variables $X_i = \sum \alpha_{ij} Y_j$ such that $k[X_1, \ldots, X_n] \cong k[Y_1,\ldots, Y_n]$ and $L_i = Y_i$ for $i \leq m$. Thus $L_1,\ldots, L_m$ is a regular sequence.  We now show that a similar fact holds for inhomogeneous linear polynomials with generic coefficients over arbitrary commutative rings, essentially because a generic matrix of size $m \times n$ with $m \leq n$ has linearly independent rows.
\begin{prop}\label{p:regular-seq-inhomog}
Let $R_0$ be a commutative ring, and let $R = R_0[\{A_{ij}\}_{1\leq i \leq m, 1\leq j \leq n}]$. Consider the linear polynomials $L_i = \sum A_{ij} X_j - c_i$, $1 \leq i \leq m$ in $R[X_1,\ldots, X_n]$, with $c_i \in R$. If $m \leq n$, then $L_1, \ldots, L_m$ is a regular sequence.
\end{prop}
\begin{proof}

We prove the statement by induction on $n$, with the base case $n=1$ being elementary.

\textbf{Claim 1.}  The sequence $L_1, \dots, L_m$ is regular over $R[A_{11}^{-1}][X_1,\ldots,X_n]$.  

\noindent To prove this, we enact the linear change of variables \[ Y_1 = L_1, \qquad Y_j = X_j \text{ for }j \geq 1.\]
 This is invertible over $R[A_{11}^{-1}]$, so  \begin{equation} \label{e:ra}
  R[A_{11}^{-1}][X_1,\ldots,X_n] = R[A_{11}^{-1}][Y_1,\ldots,Y_n]. \end{equation}
   In our new coordinates, $L_1 = Y_1$, whereas for $i > 1$ we have
  \begin{align*}
   L_i &= A_{i1}A_{11}^{-1} Y_1 + \sum_{j > 1} (A_{ij} - A_{i1}A_{11}^{-1}A_{1j}) Y_j - (c_i - A_{i1}A_{11}^{-1}c_1) \\
    &=: \Big(\sum A'_{ij} Y_j \Big) - c'_i. 
    \end{align*}
    For the first step in regularity, it is clear that $L_1 = Y_1$ is not a zero-divisor in (\ref{e:ra}).
    For the regularity of the rest of the sequence, we use the inductive hypothesis.  Let \[ S_0 = R_0[A_{i1}, A_{1j}, A_{11}^{-1}]. \]  Then recalling $Y_1 = L_1$, we see
    \begin{align}
    R[A_{11}^{-1}][Y_1,\ldots,Y_n]/(L_1) & \cong R[A_{11}^{-1}][Y_2,\ldots,Y_n]  \nonumber \\
   & =    S_0[A_{ij}]_{2 \le i, j \le n}[Y_2, \dots, Y_n]  \nonumber \\
   &  = S_0[A_{ij}']_{2 \le i, j \le n}[Y_2, \dots, Y_n].  \label{e:schange}
    \end{align}
    The last equality holds because relative to $S_0$, the change of variables $A_{ij} \mapsto A_{ij}'$ is just a translation.
By the inductive hypothesis, we know that $\bar{L}_2,\ldots, \bar{L}_m$ is a regular sequence in (\ref{e:schange}).  Claim 1 follows. 

\bigskip

\textbf{Claim 2.}  The sequence $\overline{L}_1, \overline{L}_2, \dots, \overline{L}_{m-1}$ is regular over $R/(A_{11})[X_1, \dots, X_n]$.

\noindent
Let $S_0 = R_0[A_{21},\ldots, A_{m1}][A_{m2}, \ldots ,A_{mn}][X_1]$, so that \begin{equation} \label{e:rs0}
 R/(A_{11})[X_1,\ldots,X_n] = S_0[\{A_{ij}\}_{i < m,j > 1}][X_2,\ldots,X_n]. \end{equation} Writing $c'_i = c_i - A_{i1}X_1$, we have $\bar{L}_i = \sum_{j = 2}^n A_{ij} X_j - c'_i $ in (\ref{e:rs0}). By the inductive hypothesis, $\bar{L}_1,\ldots, \bar{L}_{m-1}$ form a regular sequence in (\ref{e:rs0}). 

\bigskip

\textbf{Claim 3.}  The sequence $L_1, \dots, L_m$ is regular over $R[X_1, \dots, X_n]$.

\noindent
We induct on $m$, with the base case $m=0$ being vacuous.
We therefore suppose that $L_1,\ldots, L_{m-1}$ is a regular sequence in $R[X_1,\ldots,X_n]$. It remains to show that $L_m$ is a non-zerodivisor modulo $(L_1,\ldots, L_{m-1})$. Suppose $P L_m = \sum_{i < m} Q_i L_i$. By Claim 1, we have 
\begin{equation} \label{e:aep}
A_{11}^e P = \sum_{i < m} Q_i L_i 
\end{equation}
for some $e \geq 0$.
We need to show that there exists such an equation (with possibly different $Q_i$) with $e=0$. Therefore suppose $e > 0$.
Suppose there exists an index $i$ such that $Q_i \not \in (A_{11})$.  Let $k$ be the largest such index.
 
Since \[ \sum_{i < m} \bar{Q}_i \bar{L}_i = A_{11}^e P = 0 \quad  \text{ in } R/(A_{11})[X_1,\ldots,X_n] \]
and $\bar{Q}_i = 0$ for $i>k$, 
we have \[ \bar{Q}_k \bar{L}_k = \sum_{i < k} -\bar{Q}_i \bar{L}_i. \]
Since $k < m$, Claim 2 implies that $\bar{Q}_k \in (\bar{L}_1,\ldots, \bar{L}_{k-1})$, hence $Q_k \in (L_1,\ldots, L_{k-1}, A_{11})$. 
Write $Q_k = U_1L_1 + \cdots + U_{k-1}L_{k-1} + U_kA_{11}$, and define
\[ Q_i' = \begin{cases}
Q_i + U_i L_k & i < k \\
U_k A_{11} & i = k \\
Q_i & i > k.
\end{cases} \]
Equation (\ref{e:aep}) then becomes
\begin{equation} \label{e:aep2}
A_{11}^e P = \sum_{i < m} Q_i' L_i,
\end{equation}
where now $Q_i' \in (A_{11})$ for $i \ge k$.  Continuing the process to go from (\ref{e:aep}) to (\ref{e:aep2}), we can repeatedly change variables to decrease $k$.  In the end we obtain an equation as in (\ref{e:aep}) where  every $Q_i$ is divisible by $A_{11}$.
Since $A_{11}$ is not a zero divisor in $R[X_1, \dots, X_n]$, it can be cancelled from the equation, yielding an equation like (\ref{e:aep}) 
with $e$ replaced by $e-1$.  We can repeat this argument until we obtain an equation with $e=0$.
\end{proof}

\subsection{Application of regularity criteria to the resolution of $J'$}\label{ss:tensor-exact}

\begin{prop}\label{p:br-exact}
Let $R_0$ be a commutative ring, and let \[ R = R_0[\{ b'_i \}_{i=1}^n, \{ b_i \}_{i=1}^{n+r}, \{V_{ij}\}_{1 \leq i \leq r, 1 \leq j \leq n+r}]. \] Let 
\[  \{ 1, \ldots, n \} = S_1 \sqcup \dots \sqcup S_k \] be a partition. Define $f_i \colon R^{S_i} \lra R^2$ by $f_i(e_j) = (b_j, b'_j)$ for $j \in S_i$. Define $f_{k+1} \colon R^{r} \lra R$ by $f_{k+1}(e_i) = \sum V_{ij} b_j$. Let $I_i := \im(\det(f_i)) \subset R$, and define $I := \sum_{i=1}^{k+1} I_i$.
 Then the tensor product of complexes $\otimes_{i=1}^{k+1} R(\det(f_i))$ is a resolution of $R/I$.
\end{prop}

\begin{proof}
To see that $f_1$ is regular, $f_2 \otimes R/J_1$ is regular, all the way through $f_{k} \otimes R/(I_1 + \cdots + I_{k-1})$, we repeatedly apply Corollary~\ref{c:genericb}. The point is that the different $I_i$ use disjoint sets of variables. It remains to show that $f_{k+1}$ is regular over $R/(I_1 + \cdots + I_k)$. This is the same as showing that $L_i = \sum V_{ij} b_j$, $1 \leq i \leq r$ form a regular sequence in $R/(I_1 + \cdots + I_k)$.

Define \[ R_1 = \frac{R_0[b'_1, \ldots, b'_n, b_1,\ldots,b_n]}{(I_1 + \cdots + I_k) \cap R_0[b'_1, \ldots, b'_n, b_1,\ldots,b_n]}. \] It is easy to see that \[ R / (I_1 + \dots + I_k) = R_1[V_{ij}][b_{n+1},\ldots,b_{n+r}]. \] Set $c_i = -\sum_{j \leq n} V_{ij} b_j$. Then the linear polynomials may be written $L_i = (\sum_{j = n+1}^{n+r} V_{ij} b_j) - c_i$. Hence the $L_1, \ldots ,L_{r}$ are linear polynomials in $r$ variables with generic coefficients and non-generic constant term. By Proposition \ref{p:regular-seq-inhomog}, they form a regular sequence.
\end{proof}

We may now complete the proof of  Proposition \ref{p:resolve-j}, which as noted in \S\ref{s:outline} completes the proof of Theorem~\ref{t:comm}.  The second part of the proposition was proven in \S\ref{ss:extend-fv}, so we must only demonstrate the first part.

 Proposition~\ref{p:br-exact} implies that
the complex $C^{\bullet}$ defined in (\ref{e:cdef}) is a resolution of $R/J'$. The ring $R$ is of the appropriate form for the application of Proposition~\ref{p:br-exact} via the presentation (\ref{e:present-R}). We will briefly verify that the numerical hypothesis on the number of generic linear relations implicit in the statement of Proposition \ref{p:br-exact} is satisfied (i.e.\ that the two appearances of the variable $r$, in the number of variables $b_i$ and in the domain of $f_{k+1}$, are equal). In the matrix $D$, we have precisely one row of type III associated to each generator of the form $\boldsymbol{\rho}_{\sigma}$, $\sigma \in \mathcal{B}_v$, $v \in S \setminus \{ v_0 \}$. Therefore the number of rows of types I or II in $D$ is precisely the same as the number of generators $\boldsymbol{\rho}_i$ not associated with any $v \in S \setminus \{ v_0 \}$. This implies that the number of generators $L_i$ of $J'$ is the same as the number of variables $\boldsymbol{b}_i$ of $R$ not associated with any $v \in S \setminus \{ v_0 \}$. This is precisely the condition on the number of generic linear relations needed to apply Proposition \ref{p:br-exact} to the ideal $J'$.

\end{document}